\documentclass[review]{elsarticle}  
\usepackage{tikz}
\usepackage{amsmath,bm,algorithm,algorithmic}
\usepackage[mathscr]{euscript}
\let\euscr\mathscr \let\mathscr\relax
\usepackage{amssymb}
\usepackage{mathrsfs}
\usepackage{graphicx} 
\usepackage{natbib}
\usepackage{lineno,hyperref}
\usepackage{xcolor}
\usepackage{caption}
\usepackage{subcaption}
\usepackage{booktabs}

\begin{document}

\begin{frontmatter}

\title{Stochastic collocation approach with adaptive mesh refinement for parametric uncertainty analysis}

\author[add1]{Anindya Bhaduri}
\author[add2]{Yanyan He\footnote{Corresponding author. Email address: yanyan.he@nmt.edu}}
\author[add1]{Michael D. Shields}
\author[add1]{Lori Graham-Brady}
\author[add3]{Robert M. Kirby}
\address[add1]{Department of Civil Engineering, Johns Hopkins University, Baltimore, MD, USA}
\address[add2]{Department of Mathematics, New Mexico Institute of Mining and Technology, Socorro, NM, USA}
\address[add3]{School of Computing, University of Utah, Salt Lake City, UT, USA}




\begin{abstract}
Presence of a high-dimensional stochastic parameter space with discontinuities poses major computational challenges in analyzing and quantifying the effects of the uncertainties in a physical system. In this paper, we propose a stochastic collocation method with adaptive mesh refinement (SCAMR) to deal with high dimensional stochastic systems with discontinuities. Specifically, the proposed approach uses generalized polynomial chaos (gPC) expansion with Legendre polynomial basis and solves for the gPC coefficients using the least squares method. It also implements an adaptive mesh (element) refinement strategy which checks for abrupt variations in the output based on the second order gPC approximation error to track discontinuities or non-smoothness. In addition, the proposed method involves a criterion for checking possible dimensionality reduction and consequently, the decomposition of the full-dimensional problem to a number of lower-dimensional subproblems. Specifically, this criterion checks all the existing interactions between input dimensions of a specific problem based on the high-dimensional model representation (HDMR) method, and therefore automatically provides the subproblems which only involve interacting dimensions. The efficiency of the approach is demonstrated using both smooth and non-smooth function examples with input dimensions up to 300, and the approach is compared against other existing algorithms.
\end{abstract}

\begin{keyword}
Generalized polynomial chaos, stochastic collocation, adaptive mesh refinement, interaction check 
\end{keyword}

\end{frontmatter}

\section{Introduction}
Computer-based simulations are widely used for predicting the behavior of physical systems. However, due to uncertainties in the system and the simulation process, such as the inherently stochastic nature of some system parameters, boundary conditions or excitations and a lack of understanding of the true physics, predictions inevitably deviate from reality. Therefore, understanding and quantifying the uncertainty in simulations is necessary in order to incorporate potential variability into these predictions. 

One of the main aspects of uncertainty quantification (UQ) is uncertainty propagation, also called forward UQ. It aims to quantify uncertainty in the model outputs that results from uncertainty in the model inputs, which are usually represented using random variables with an associated probability distribution. The goal is therefore to estimate the response surface, probability density function (PDF) or statistical moments for the model outputs efficiently. Probabilistic approaches have been relatively well-developed for forward UQ. For example, the most popular technique is the Monte Carlo method, which is robust, simple to understand, easy to implement, and typically serves as a baseline against which other methods are compared. However, it may require a large number of model evaluations to reach the desired accuracy due to its slow convergence rate. 

Other efficient methods have been proposed to achieve a higher convergence rate and consequently reduce the computational cost. Polynomial chaos (PC) expansion is one such method which represents the output of interest by the expansion of orthogonal polynomials (with respect to positive weight measure) in the stochastic input space. It is based on the homogeneous chaos theory by Wiener \cite{wiener1938homogeneous} where a Gaussian process was essentially expressed by a set of Hermite polynomials. Ghanem and Spanos \cite{ghanem2003stochastic} have coupled this approach with finite element methods to effectively model uncertainty in solid mechanics problems. The generalized polynomial chaos (gPC) \cite{xiu2002wiener,xiu2003modeling} method makes use of different types of orthogonal polynomials in the Askey scheme \cite{askey1985some} as the bases to approximate random functions/processes. It is capable of reaching fast convergence for smooth functions when the PDF of the random variables is identical to the weighting function of the orthogonal polynomials from the Askey scheme. This idea has been further extended to arbitrary random distributions \cite{witteveen2006modeling,wan2006beyond}. The gPC coefficients in the above works are determined by performing Galerkin projection on the model equations. Its intrusive nature requires the modification of the deterministic simulation code, which could be a difficult and time-consuming task. 

By contrast, non-intrusive methods use the deterministic simulation code directly without requiring any modifications, which makes them more applicable to complex systems. For example, Xiu \cite{xiu2007efficient} proposed a gPC scheme based on the stochastic collocation method, where the gPC coefficients are obtained using the discrete projection approach. Babuska et.al. \cite{babuvska2007stochastic} used Gauss quadrature points to sample low dimensional random spaces and perform tensor product interpolation using 1-D basis functions. Tensor grid approaches suffer from the so-called `curse of dimensionality' \cite{cools2002advances} as there is an exponential rise in the required number of full model evaluations with the increase in dimensionality of the input space. To alleviate this problem to some extent, sparse grid \cite{smolyak1963quadrature,bungartz2004sparse} based interpolations \cite{xiu2005high,nobile2008sparse} have been performed with the global Lagrange polynomial basis as the interpolant in the random space. However, these global approaches may not be suitable for tracking local steepness or discontinuities in the random space, and the approximation may fail to converge to the true value.

To deal with non-smooth functions, multi-element schemes have been proposed for both intrusive and non-intrusive methods. Wan and Karniadakis \cite{wan2005adaptive} developed a multi-element generalized polynomial chaos (MEgPC) scheme based on the stochastic Galerkin method to handle the issue of discontinuities in the output response and long-term integration of stochastic differential equations. This approach adaptively splits the actual input domain into smaller subdomains by calculating the relative error in variance along each dimension and maintaining a relatively low polynomial order (less than 10) in critical subdomains. However, as an intrusive approach, it requires modification of the deterministic simulation code. Foo et. al. \cite{foo2008multi} introduced the non-intrusive multi-element probabilistic collocation method (MEPCM) with Lagrange polynomial basis to efficiently treat problems characterized by strong non-linearities or discontinuities and long-term integration. The criterion for adaptively splitting the input domain is similar to that in the MEgPC scheme. 

Both the Galerkin and collocation versions of the multi-element gPC scheme are still dimension-dependent, since both the number of subdomains and the number of terms in the gPC expansion increase rapidly with the increase in dimensionality of the stochastic input. To mitigate the issue of high computational cost associated with the element decomposition in high dimensional problems, Foo and Karniadakis \cite{foo2010multi} developed the MEPCM-A method, which combines the MEPCM with the high dimensional model representation (HDMR) \cite{sobol2003theorems}. The HDMR represents a function as a hierarchical additive combination of lower dimensional functions starting from a one-dimensional input space to a full-dimensional input space. A way to estimate the correlation functions is to use the cut-HDMR approach \cite{rabitz1999efficient}. In the MEPCM-A approach, a high-dimensional stochastic problem is reduced to a series of low-dimensional problems by truncating the terms in the HDMR up to a certain dimensionality, $\nu$, followed by the application of the MEPCM approach to each of these subproblems with maximum dimensionality $\nu$. Parameter $\nu$ is generally chosen to be small enough compared to the high dimensionality of the original problem that element decomposition is not computationally prohibitive. Another important parameter in the MEPCM implementation is the number of points, $\mu$, in the interpolation rule. Parameters $\nu$ and $\mu$ are pre-fixed without regard to the actual order of interaction among the input parameters. For problems with high nominal dimensions but low effective dimensions (i.e. only a few input variables strongly influence the response), the method proves to be efficient. However, the choice of a proper value for $\nu$ of the subproblems needs more exploration. In addition, once $\nu$ is prescribed, all the interaction terms up to order $\nu$ in the HDMR are considered. Consequently, for complex systems with strong input interactions, $\nu$ may be chosen to be large for satisfactory error estimates and thus the number as well as the dimensionality of the subproblems could become prohibitively large. Even with a small value of $\nu$, the number of interaction terms can become very large for very high dimensional problems. Moreover, the model output may not be sensitive to some interaction terms with order upto $\nu$, and thus a significant number of unnecessary sub-problems are considered which increases the computational cost.

Approaches \cite{le2004uncertainty,mathelin2003stochastic,babuvska2007stochastic} based on local bases have also been proposed to deal with non-smoothness in the random space. Klimke and Wohlmuth \cite{klimke2005algorithm} developed a sparse grid collocation interpolation scheme based on piecewise linear basis functions, which has the ability to resolve discontinuities in the response surface but suffers from slow convergence rates because of global refinement of the sparse grid. The approach is based on hierarchical sparse grid points where points are added in successive depth levels. The error indicator is known as the hierarchical surplus and acts as a stopping criterion for the algorithm. Ma and Zabaras \cite{ma2009adaptive} used a similar approach called adaptive sparse grid collocation (ASGC) but also incorporated an adaptive strategy that enables a local sparse grid refinement around the discontinuity region, which helps enhance the convergence rate. The ASGC approach checks the hierarchical surplus values at each point in the current depth level and creates new points in the next depth level only in the neighborhood of points whose surplus error exceeds the tolerance value. 
The approach is restricted to uniform grid points because of the adaptivity criterion. For the purpose of tracking discontinuities, ASGC uses piecewise linear basis function. This may lead to a slow convergence for the regions where the approximating response surface are smooth. To tackle high dimensional stochastic problems, Ma and Zabaras \cite{ma2010adaptive} combined a dimension-adaptive version of HDMR with ASGC (HDMR-ASGC). Initially, the importance of the component functions in HDMR are estimated through a weight measure which is expressed as the integral value of a component function of certain order with respect to the sum of the integral values of all lower order component functions. Component functions with weight measures higher than a predefined error threshold are the ones considered important. ASGC is then applied to each of the lower dimensional sub-problems corresponding to the important component functions. The error indicator used in HDMR-ASGC is a function of the integral value of the basis function as well as the hierarchical surplus. It is different from the original ASGC approach \cite{ma2009adaptive} which uses only the surplus value as the error indicator.

In this paper, we propose a method of stochastic collocation with adaptive mesh refinement (SCAMR). Specifically, the proposed approach uses generalized polynomial chaos (gPC) expansion with Legendre polynomial basis and solves for the gPC coefficients using the least squares method. 
It also implements an adaptive mesh (element) refinement strategy to track any discontinuities or non-smoothness in the output. The adaptive criteria associated with the mesh refinement strategy check for abrupt variations in the output based on the observed error from a second order gPC approximation. SCAMR further introduces a criterion for possible dimensionality reduction, allowing for decomposition of the full-dimensional problem to a number of lower-dimensional subproblems. This criterion checks all the existing interactions between input dimensions of a specific problem based on HDMR, and consequently provides the subproblems which only involve interacting dimensions.

The paper is organized as follows: Section 2 presents the general framework for a stochastic problem. In Section 3, we discuss the proposed method of stochastic collocation with adaptive mesh refinement in detail. In Section 4, we demonstrate the effectiveness and efficiency of the proposed approach using various numerical examples compared to the ASGC, the HDMR-ASGC as well as the MEPCM-A approach. We finally conclude the paper with a discussion in Section 5. 

\section{Problem Definition}
Let the triplet ($\Omega,\mathcal{F},\mathcal{P}$) represent a complete probability space, where $\Omega$ corresponds to the sample space of outcomes, $\mathcal{F} \subset 2^{\Omega}$ is the $\sigma$-algebra of measurable events in $\Omega$, and $\mathcal{P} : \mathcal{F} \rightarrow [0,1] $ is the probability measure. Let $\bm{\xi} =\{\xi_1(\omega),  \xi_2(\omega), \dots, \xi_{n}(\omega)\}: \Omega \rightarrow \Xi \in \mathbb{R}^n$ be a set of $n$ independent random variables, which characterize the uncertainty in the system. In the current work, we assume that the random variables $\bm{\xi}_i$ follow uniform distribution with a constant PDF $p(\bm{\xi})=\rho_{\bm{\xi}}; \ \bm{\xi} \in [a_1,b_1] \times [a_2,b_2] \times .... \times [a_n,b_n]$. Let $ \bm{x} \in D \subset \mathcal{R}^d \ (d \in \{1,2,3\})$ be the spatial variable, and $t \in (0,T]$ ($T>0$) be the temporal variable.

Consider a general partial differential equation
\begin{equation}\label{govern}
\begin{cases}
u_{t}(x,t,\bm{\xi})=\mathcal{L}(u;\bm{x},t,\bm{\xi}), & D\times (0,T] \times \Xi\\
 \mathcal{B}(u; \bm{x},t,\bm{\xi})=0, & \partial{D}\times[0,T]\times \Xi,\\
u=u_{0}, & \bar{D}\times\{t=0\} \times \Xi,
\end{cases}
\end{equation}
where $\mathcal{B}$ is the operator for the boundary conditions, $\mathcal{L}$ is the differential operator, $D$ is the spatial domain, and $u=u_{0}$ is the initial condition.  The problem is assumed to be well-posed in parameter space $ \Xi$.
The model output $u(\bm{x},t,\bm{\xi})$ is the quantity of our interest. For the convenience of notation, we do not consider the dependence of solution on the spatial and time variables $\bm{x}$ and $t$, and only discuss the problem for any fixed $\bm{x} \in D$ and $t \in (0,T]$. As mentioned in \cite{chen2015efficient}, this is standard in the UQ literature. Our goal is to quantify the uncertainty in the quantity of interest $u(\cdot, \bm{\xi}): \Xi \to \mathbb{R}$, due to the uncertainty in the input variables $\bm{\xi}$. Without loss of generality, we consider scalar model output. 

\section{Stochastic Collocation with Adaptive Mesh Refinement}

In this section, we propose a stochastic collocation method with adaptive mesh refinement (SCAMR). Specifically, SCAMR adopts a mesh refinement scheme with a proposed criteria that checks for discontinuities or abrupt variations in the response surface, as well as interactions between different input dimensions. Details are provided in the following subsections.

\subsection{Generalized Polynomial Chaos Based Stochastic Collocation}
Let $u(\bm{\xi}) \in L_2(\Xi)$ be a square-integrable function of the $n$-dimensional random vector $\bm{\xi}$ which can be represented using the generalized polynomial chaos expansion as
\begin{equation}
u(\bm{\xi}(\omega))=\sum_{i=0}^{\infty} \hat{u}_i\Phi_i(\bm{\xi}(\omega)),
\end{equation}
where $\hat{u}_i$ are the gPC coefficients and $\Phi_i$ are the Legendre polynomials for uniform $\bm{\xi}$ \cite{xiu2002wiener}.

For numerical calculations, the series is truncated to $N+1$ terms to approximate the exact output $u(\bm{\xi}(\omega))$ with polynomial order $p$ 
\begin{equation}\label{gPC_trunc}
u_p(\bm{\xi}(\omega))=\sum_{i=0}^{N} \hat{u}_i\Phi_i(\bm{\xi}(\omega)), \qquad N+1=\frac{(n+p)!}{n!p!},
\end{equation}
where
\begin{equation}
\hat{u}_i=\frac{1}{E[\Phi_i^2]} \int_\Xi u(\bm{\xi}) \Phi_i(\bm{\xi})  \rho(\bm{\xi}) d\bm{\xi}.
\end{equation}
With collocation methods, the gPC coefficients $\hat{u}_i$ can be obtained using discrete projection as
\begin{equation}\label{discrete-proj}
\hat{u}_i=\frac{1}{E[\Phi_i^2]}\sum_{j=1}^{M}u(\bm{\xi}^j)\Phi_i(\bm{\xi}^j)\alpha^j, \qquad i=0,1,\dots,N,
\end{equation}
where $\{\bm{\xi}^j,\alpha^j\}_{j=1}^M$ are sets of quadrature points and their corresponding weights. 

Another collocation method for estimating the gPC coefficients utilizes interpolation on the pairs $\{\bm{\xi}^j, u(\bm{\xi}^j)\}_{j=1}^{N+1}$. The gPC coefficient vector $\hat{\textbf{u}} =\{\hat{u}_0, \hdots, \hat{u}_N\}$ is estimated by solving the following linear system
\begin{equation}
\sum_{i=0}^{N} \hat{u}_i\Phi_i(\bm{\xi}^j)=u(\bm{\xi}^j), \forall j=1,2,\dots,N+1.\nonumber
\end{equation}
The interpolation method may not produce a proper approximation if $u(\bm{\xi}^j)$ is corrupted by observational or measurement errors. The projection method, on the other hand, produces the best approximation in the weighted $L_2$ norm \cite{xiu2010numerical}. However, the quadrature nodes used in the discrete projection method have restrictions, such as the structure of the nodes and the number of the nodes. 

To allow more flexibility, in terms of the location and the number of nodes, we estimate the vector of gPC coefficients by solving the following least squares problem using $M \ ( \ > N+1)$ sets of points:
\begin{equation}\label{gPC_leastsq}
\hat{\textbf{u}}=\arg\min_{\tilde{u}} \| \sum_{i=0}^{N} \tilde{u}_i\Phi_i(\bm{\xi})-u(\bm{\xi})\|_2
\end{equation}
where $\tilde{\textbf{u}}=\{\tilde{u}_0,\tilde{u}_1,\hdots,\tilde{u}_N\}$ is an arbitrary gPC coefficient vector which converges to the desired vector $\hat{\textbf{u}}=\{\hat{u}_0,\hat{u}_1,\hdots,\hat{u}_N\}$ through the minimization in Eq. (\ref{gPC_leastsq}).
Consequently, the approximated output $u_p$ is estimated using Eq. (\ref{gPC_trunc}). It is to be noted here that the set of $M$ points may have an unstructured arrangement in the input space.

\subsection{Decomposition of Random Space}
In this section, we introduce the standard decomposition method for random input space, where the $L_2$ error of the global approximation has been proven to be bounded by the local $L_2$ error approximations in the elements \cite{wan2005adaptive}. We assume a hypercube input domain in our present work. Without the loss of generality, we consider the original stochastic space as $\Xi=[-1, 1]^n$. It is then decomposed into $n_e$ non-overlapping and space-filling elements $\Xi_k$: $\cup_{k=1}^{n_e} \Xi_k = \Xi$, $\Xi_m \cap \Xi_k = \emptyset$ for $m\neq k$ and $m,k \in [1,2,\dots, n_e] $. If $a_i^k$ and $b_i^k$ denote the minimum and maximum bounds of element $\Xi_k$ along dimension $i$ ($1\leq i \leq n$), $\Xi_k$ is the tensor product given by
\begin{equation}
\Xi_k=[a_1^k,b_1^k)\times[a_2^k,b_2^k)\times..........\times[a_n^k,b_n^k).
\end{equation}

Let the local input random vector in each element be defined as $\bm{\xi}^k=[\xi_1^k,\xi_2^k,\dots,\xi_n^k]$. For the purpose of applying the gPC formulation on each element locally, the local random vector can be transformed to a new random vector $\bm{\eta} \in [-1, 1]^n$  such that  $\bm{\eta}=F_k(\bm{\xi}^k)=[\eta_1,\eta_2, \dots, \eta_n]$. The transformation is a simple scaling relationship between the $[-1,1]^n$ domain and the particular $\Xi_k$ domain:
\begin{equation}
F_k:\eta_i=-1+ \frac{2}{b_i^k-a_i^k}(\xi_i^k-a_i^k), \ \ \forall i=1,2,...,n
\end{equation}

\subsection{Adaptive Criteria}
The SCAMR algorithm uses adaptive approaches for two purposes: detection of abrupt variations in the output function for non-smoothness and reduction of the high-dimensional input parameter space to a subset of interacting dimensions. Each of these are described in the following subsections.
\subsubsection{Criterion for Detecting Abrupt Variation in One Dimension}
In the current work, we propose to use first or second order Legendre polynomials to efficiently approximate any general response function with local abruptness or discontinuities. In any domain where the function deviates significantly from a second order polynomial approximation, we decompose the domain further. Specifically, we consider the output variation along the centerline (straight line passing through the center of the domain) along each dimension one at a time with the rest of the dimensions fixed at their midpoints. For example, let $\Gamma$ be a given $n$-dimensional domain (element) such that $\Gamma=[a_1,b_1)\times[a_2,b_2)\times.....\times[a_n,b_n)$. For the $i$-th dimension, let $\bm{z}=\{z_1,\hdots, z_m\}$ be $m$ Chebyshev points of depth level $l$ in the range $[a_i,b_i)$ such that $m=2^{l}+1$. In this study, depth level $l=2$ is taken and hence $m=5$. Then the set of input points along the centerline in the $i$-th dimension is $\bm{\xi}^{(i)}=\{\bm{\xi}^{(i)}_{1},\bm{\xi}^{(i)}_{2},\hdots,\bm{\xi}^{(i)}_{m}\}$, where each $n$-dimensional point is $\bm{\xi}^{(i)}_{j}=\{\frac{a_1+b_1}{2}, \hdots, \frac{a_{i-1}+b_{i-1}}{2}, z_j, \frac{a_{i+1}+b_{i+1}}{2}, \hdots, \frac{a_n+b_n}{2}\},\forall j \in \{1,2,\dots,m \}$. Let $\bm{u}^{(i)}=\{u^{(i)}_{1},u^{(i)}_{2},\hdots,u^{(i)}_{m}\}$ be the corresponding set of $m$ exact outputs and $\bm{u}^{(i)}_p=\{u^{(i)}_{p,1},u^{(i)}_{p,2},\hdots, u^{(i)}_{p,m}\}$ be the corresponding $1$-D second-order  gPC approximation along the $i$-th dimension for the current domain. The model output can then be reasonably approximated as quadratic if
\begin{equation}\label{1D-check}
\|\bm{u}^{(i)}_p-\bm{u}^{(i)}\|_{\infty} < \epsilon_1,
\end{equation}
where $\epsilon_1$ is an error tolerance parameter. If criterion (\ref{1D-check}) is not satisfied, the $i$-th dimension is considered critical. All the critical dimensions are then stored in descending order of the error magnitude obtained from criterion (\ref{1D-check}) and the domain is further decomposed along the center of the two most critical dimensions. The domain subdivision is repeated for every newly formed element until the stopping criteria are satisfied.
\subsubsection{Criterion for Dimensionality Reduction}\label{cr: interaction}
The second criterion helps in achieving dimensionality reduction. It decomposes the original full-dimensional problem to a number of lower dimensional problems by identifying the absence of interactions between input dimensions with respect to the output of interest. This criterion is checked at two levels and takes advantage of the significant gains in computational efficiency by dealing with low-dimensional functions.

{\it{\textbf{First level criterion}}}. At the first level, a dimension $i$ is assumed non-interacting with others if 
\begin{equation}\label{no-variation}
||\bm{u}^{(i)}-u_c||_{\infty} < \epsilon_1,
\end{equation}
where $\bm{u}^{(i)}$ is the centerline output vector along the $i$-th dimension (introduced earlier) and $u_c$ is the exact output value at the center point of the input domain $\Xi$. 
By implementing this first level criterion, the full-dimensional problem will be decomposed to a $r\leq n$ dimensional and $n-r$ one-dimensional problems, where the one-dimensional problems depend on the input random variables which do not interact with others.

{\it{\textbf{Second level criterion}}}. At the second level, we further decompose the $r$-dimensional problem to a number of lower-dimensional sub-problems by verifying ${r}\choose{2}$ pairwise interactions in the $r$-dimensional domain. All higher dimensional interactions between the input dimensions are derived from the pairwise interaction results. This second level criterion is derived from the HDMR representation \cite{rabitz1999efficient,rabitz1999general} and the details are provided in the following.\\
{\it{Pairwise non-interaction criterion derivation}}.  
Let $f(\bm{Y})=f(Y_1,Y_2,....,Y_n)$ be an $n$-dimensional function. Following the notation in \cite{ma2010adaptive}, the general expression of the High Dimensional Model Representation (HDMR) for the function is given by
\begin{align}\label{HDMReq}
f(\bm{Y})&=f_0+\sum_{i=1}^nf_{i}(Y_{i})+\sum_{1 \le i_1<i_2 \le n}f_{i_1i_2}(Y_{i_1},Y_{i_2})+\hdots \nonumber \\
& +\sum_{1 \le i_1<..i_s \le n}f_{i_1....i_s}(Y_{i_1},...,Y_{i_s})
+...... +f_{12...n}(Y_{1},Y_{2},...Y_{n})
\end{align}
where $f_0$ is a constant zeroth order function, $f_i()$ denotes a one-dimensional function, $f_{i_1i_2}()$ is a two-dimensional function and so on.

As seen from Eq. (\ref{HDMReq}), the HDMR breaks down the function $f(\bm{Y})$ into individual contributions from all possible orders of interactions among the dimensions. For example, $f_i(Y_i)$ represents how input $Y_i$ influences $f(\bm{Y})$ keeping the other input dimensions fixed. The third term $f_{i_1i_2}(Y_{i_1},Y_{i_2})$ represents the combined contribution of inputs $Y_{i_1}$ and $Y_{i_2}$ towards $f(\bm{Y})$ after their individual contributions have been accounted for through $f_i(Y_i)$. All dimensions except $Y_{i_1}$ and $Y_{i_2}$ are kept fixed in this case. Similarly, $f_{12...n}(Y_{1},Y_{2},...Y_{n})$ denotes the contribution of all inputs taken together towards $f(\bm{Y})$ after having accounted for all lower dimensional function contributions.

Cut-HDMR \cite{li2001high1,li2001high2} is an efficient technique for estimating the component functions in $f(\bm{Y})$ which involves evaluating $f(\bm{Y})$ on lines, planes and hyper-planes (or cuts) passing through a “cut” center $\bm{c}$ which is a point in the input variable space. The choice of $\bm{c}$ is important as it influences the convergence of the HDMR expansion. It has been shown \cite{xu2004generalized} that a suitable choice of $\bm{c}$ can be the mean of the input random vector.

The component functions \cite{shan2010metamodeling} are given by:
\begin{align}
f_0=f(\bm{c}) 
\end{align}
\begin{align}\label{1st order interaction}
f_i(Y_i)=f(Y_i,\bm{c}^{\{i\}})-f_0 \ \ \ \forall i \in \{1,2,\hdots,n\}
\end{align}
\begin{align}\label{2nd order interaction}
f_{i_1i_2}(Y_{i_1},Y_{i_2})=f(Y_{i_1},Y_{i_2},\bm{c}^{\{i_1,i_2\}})-f_{i_1}(Y_{i_1})-f_{i_2}(Y_{i_2})-f_0,
\end{align}
$\ \ \ \ \ \ \ \ \ \ \ \ \ \ \ \ \ \ \ \ \ \ \ \ \ \ \ \ \ \ \ \ \ \ \ \ \ \ \ \ \ \ \ \ \ \ \forall i_1, i_2 \in \{1,2,\hdots,n\},$ such that $i_1< i_2$
\begin{align}
f_{i_1i_2i_3}(Y_{i_1},Y_{i_2},Y_{i_3})&=f(Y_{i_1},Y_{i_2},Y_{i_3},\bm{c}^{\{i_1,i_2,i_3\}})-f_{i_1i_2}(Y_{i_1},Y_{i_2})
-f_{i_1i_3}(Y_{i_1},Y_{i_3}) \nonumber \\
&\qquad{}
-f_{i_2i_3}(Y_{i_2},Y_{i_3})-f_{i_1}(Y_{i_1})-f_{i_2}(Y_{i_2})-f_{i_3}(Y_{i_3})-f_0,
\end{align}
$\ \ \ \ \ \ \ \ \ \ \ \ \ \ \ \ \ \ \ \ \ \ \ \ \ \ \ \ \ \ \ \ \ \ \ \ \forall i_1, i_2, i_3 \in \{1,2,\hdots,n\},$ such that $i_1< i_2 < i_3$\\
$\ \ \ \ \ \ \ \ \ \ \ \ \ \ \ \ \ \ \ \ \vdots$
\begin{align}
f_{12 \hdots n}(Y_{1},Y_{2},\hdots,Y_{n})&=f(\bm{Y})-f_0-\sum_{i=1}^nf_{i}(Y_{i_1})-\sum_{1 \le i_1<i_2 \le n}f_{i_1i_2}(Y_{i_1},Y_{i_2})\nonumber \\
&\qquad{}
-\hdots-\sum_{1 \le i_1<..i_{n-1} \le n}f_{i_1 \hdots i_{n-1}}(Y_{i_1},\hdots, Y_{i_{n-1}})
\end{align}
where $\bm{c}^{\{i\}}=\bm{c}$\textbackslash $\{Y_i\},\bm{c}^{\{i_1,i_2\}}=\bm{c}$\textbackslash $\{Y_{i_1},Y_{i_2}\},\bm{c}^{\{i_1,i_2,i_3\}}=\bm{c}$\textbackslash $\{Y_{i_1},Y_{i_2},Y_{i_3}\}$. For sets $A$ and $B$, $A$\textbackslash$B$ denotes a set with only those elements in $A$ that are not included in $B$.

Using the HDMR representation, we will now derive the non-interaction criterion for dimensionality reduction. 
In the proposed method, we consider only pairwise interactions of inputs. We thus concentrate on the second order (2-dimensional) component function given by Eq. (\ref{2nd order interaction}). 
Combining Eq. (\ref{1st order interaction}) with Eq. (\ref{2nd order interaction}), we can write,
\begin{align}
f_{i_1i_2}(Y_{i_1},Y_{i_2})=f(Y_{i_1},Y_{i_2},\bm{c}^{\{i_1i_2\}})-f(Y_{i_1},\bm{c}^{\{i_1\}})-f(Y_{i_2},\bm{c}^{\{i_2\}})+f_0,
\end{align}
For a given error tolerance $\epsilon_2$, dimensions $i_1$ and $i_2$ can be considered non-interacting if the second order component function $f_{i_1i_2}(Y_{i_1},Y_{i_2})$ is considered negligible, i.e.,
$f_{i_1i_2}(Y_{i_1},Y_{i_2}) \le \epsilon_2$. This implies,
\begin{align}\label{non-interaction}
f(Y_{i_1},Y_{i_2},\bm{c}^{\{i_1,i_2\}})-f(Y_{i_1},\bm{c}^{\{i_1\}})-f(Y_{i_2},\bm{c}^{\{i_2\}})+f_0 \le \epsilon_2. 
\end{align}
Eq. (\ref{non-interaction}) is the pairwise non-interaction criterion.\\
\begin{figure}[h]
\centering
\includegraphics[width=0.7\linewidth, height=6cm]{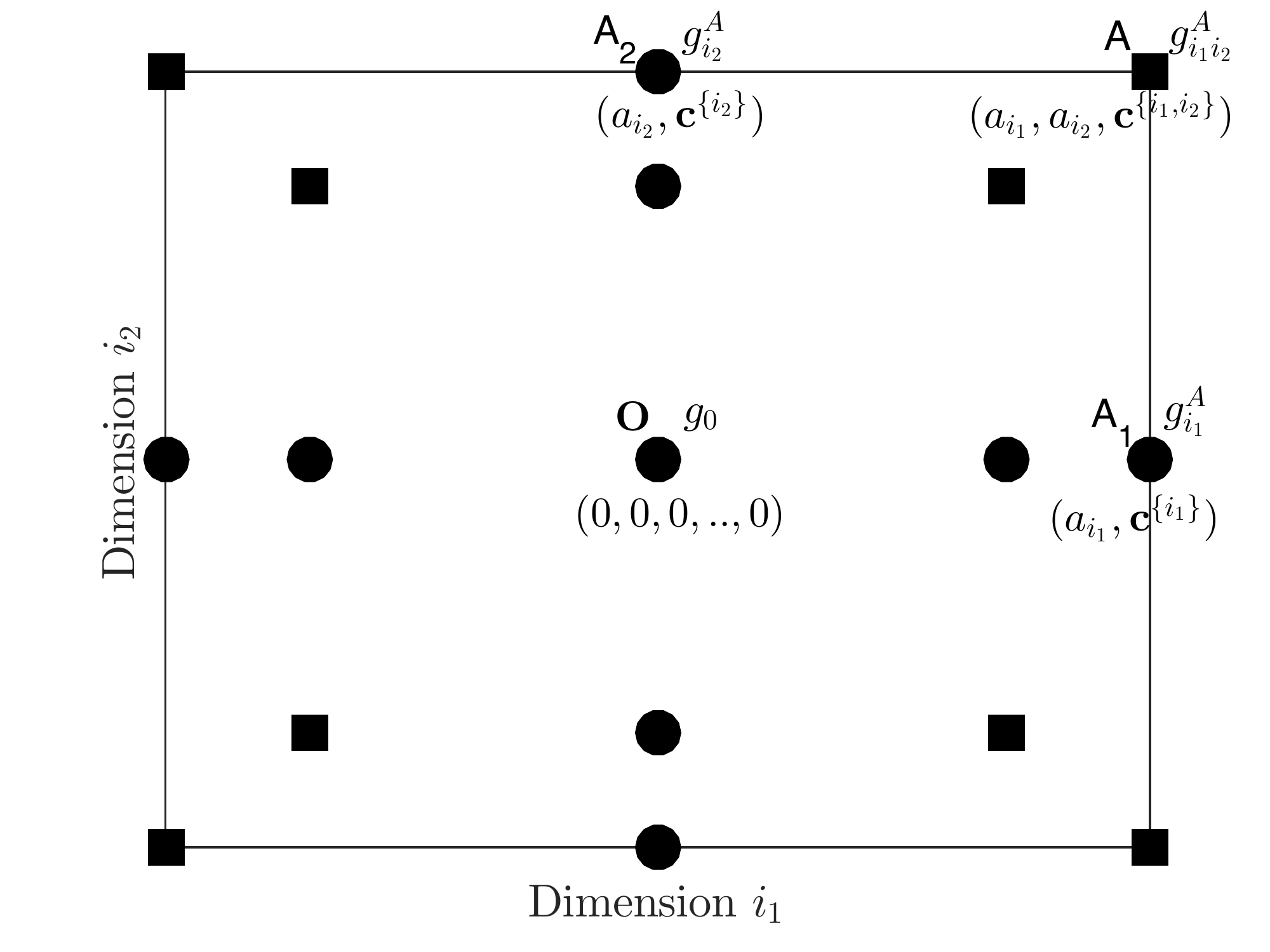} 
\caption{Square points denote new points introduced for the interaction check between dimensions}
\label{fig:subim0131}
\end{figure}
\indent Let us take a two-dimensional input domain as an example (see Fig. \ref{fig:subim0131}), where the input domain is projected from a higher $n$-dimensional input space with all the dimensions fixed at the mean of their respective ranges 
except those two dimensions ($i_1$ and $i_2$). The cut center is given by $c=\{0,0, \dots ,0\}$ and is denoted by point $O$ in Fig. \ref{fig:subim0131}. All the square points in Fig. \ref{fig:subim0131} are used to test for interaction between the two dimensions. The exact values at all the square points are calculated by full model evaluations and compared with the values at those points obtained assuming both dimensions are non-interacting. For example, assume the exact value at point A is $g_{i_1i_2}^A$ and the approximated value at A assuming non-interaction is given by
$g_{i_1i_2}^{approx,A}=g_{i_1}^A+g_{i_2}^A-g_0$. The output values $g_{i_1}^A$ and $g_{i_2}^A$ 
correspond to input points at $A_1$ and $A_2$ which are orthogonal projections of $A$ on axes $i_1$ and $i_2$ respectively passing through 
point $O$  and $g_0$ is the corresponding output value. Let $\textbf{g}^{\ true}_{i_1i_2}$ be the true output vector corresponding to the square points and $\textbf{g}^{\ approx}_{i_1i_2}$ be the corresponding approximate output vector obtained from the outputs at the circular points such that $\textbf{g}^{\ approx}_{i_1i_2}=\textbf{g}_{i_1}+\textbf{g}_{i_2}-g_0$. 
Then, Eq. (\ref{non-interaction}) is considered satisfactory if
\begin{equation}\label{non-interac}
\|\textbf{g}^{true}_{i_1i_2}-\textbf{g}^{approx}_{i_1i_2}\|_{\infty} \leq \epsilon_2,
\end{equation}

As mentioned earlier, using the knowledge about each of the pairwise (2-dimensional) interactions, we derive all the possible higher dimensional interactions. For example, we consider a 5-dimensional stochastic function where $\{Y_1,Y_2,\hdots,Y_5\}$ are the input dimensions. If only pairs $\{Y_1,Y_2\},\{Y_2,Y_3\}$ and $\{Y_3,Y_1\}$ out of total ${5}\choose{2}$ $=10$ pairs are interacting based on the criterion Eq. (\ref{non-interaction}), we decompose the full five-dimensional problem into a three-dimensional problem in the space of $\{Y_1,Y_2,Y_3\}$ and two one-dimensional problems in the space of $Y_4$ and $Y_5$, respectively.

{\it{\textbf{Sub-dimensional representation}}}. After checking criteria in Eqs. (\ref{no-variation}) and (\ref{non-interaction}), an $n$-dimensional problem can be potentially reduced to a set of lower dimensional problems as mentioned in the beginning of this section. We discuss next the effects of applying the two criteria in successive steps and how to represent the full-dimensional function in terms of a number of lower dimensional functions. At first, using criterion (\ref{no-variation}), an $n$-dimensional input domain $\Xi$ of dimension index set $D=\{1,2,...,n\}$ can be potentially reduced to a group of $N_R$ non-interacting lower dimensional input domains of dimension index set $R=\{R_1,R_2,R_3, ...., R_{N_R} \}$ with 
$|R_1| =r $, $|R_j|=1 \ (\forall j=2,....,N_R)$, $\cup_{i=1}^{N_R}{R_i}=D$, $\cap_{i=1}^{N_R}{R_i}=\O$ and $N_R=n-r+1$.

In the next step, using criterion (\ref{non-interaction}), the $R_1$ sub-dimensional problem can be further reduced to a group of $N_Q$ lower dimensional input domains with dimension index set $Q=\{Q_1,Q_2,......,Q_{N_Q}\}$ such that $\cup_{i=1}^{N_Q}{Q_i}=R_1$. Thus, in total, an $n$-dimensional problem can be reduced to $N_{S}$ lower dimensional input domains with dimension index set $S=\{R_2,R_3,.....,R_{N_R},Q_1,Q_2,....,Q_{N_Q}\}=\{S_1,S_2,.....,S_{N_S}\}$ such that $N_S=N_R+N_Q-1$. The $N_S$ index sets can be overlapping such that $\cap_{i=1}^{N_S} S_i \ne \O$. 
In case of overlapping, common dimension indices will be present among different elements in $S$. These common dimension indices form $N_T$ additional low dimensional domains of dimension set $T=\{T_1,T_2,...,T_{N_T}\}$, where $T=\{S_i \cap S_l\}$ \textbackslash \O $, \ \ \forall i,l \in \{1,....,N_S\}$ such that $i < l$. These additional low dimensional functions can be called ``corrective" dimension index sets introduced in order to account for the overlapping in $S$. Each of the ``corrective" sets has an associated constant factor $U_j \ (\forall j=1,2,...,N_T)$, which equals the difference between frequency of its occurrence in $S$ and the frequency of its occurrence in $T$. The frequency of occurrence of an index set in $S$ or $T$ is the number of times an index set features in $S$ or $T$ by itself or as a subset in a larger index set. There is also a constant factor $V$ associated with $f_0$, the function value at the cut center. In case of no overlapping of elements in $S$, i.e., $\cap_{i=1}^{N_S} S_i = \O$, then $T=\{\O\}$ and $N_T=0$. The function can thus have an HDMR-like representation and is given by
\begin{align}
f(Y_1,Y_2,....Y_n)&=\sum_{i=1}^{N_S}h_i(\bm{Y}_{S_i},\bm{c}^{S_i})-\sum_{j=1}^{N_T}U_jp_j(\bm{Y}_{T_j},\bm{c}^{T_j})
-Vf_0,
\end{align}
where $\bm{Y}_{S_i}$ is the set of input variables with the elements in $S_i$ as the indices, $\bm{Y}_{T_j}$ is the set of input variables with the elements in $T_j$ as the indices, $h_i()$ is an $|S_i|$-dimensional function, $p_j()$ is a $|T_j|$-dimensional function, $U_j$ and $V$ are integer constants where $V=N_S-\sum_{j=1}^{N_T} U_j-1$.

As an example, consider an $8$-dimensional function $f(\bm{Y})$. 
It is assumed that from criterion (\ref{no-variation}), each of the last $r=3$ dimensions is identified to be non-interacting with the remaining $(n-1)=7$ dimensions. We thus have the following set of non-interacting group of dimensions:
\begin{equation}
R=\{\{1,2,\hdots,5\},\{6\},\{7\},\{8\}\}, \nonumber 
\end{equation}
and the function can now be described by:
\begin{align}\label{Eq: example-no-var}
f(\bm{Y})&=f(Y_1,Y_2,\hdots,Y_{8}) \nonumber \\
&=g_{0}(Y_1,Y_2,\hdots,Y_{5},\bm{c}^{\{1,2,\hdots,5\}})+h_{1}(Y_{6},\bm{c}^{\{6\}})\nonumber\\
\qquad{}
& \ \ \ \ \ +h_{2}(Y_{7},\bm{c}^{\{7\}})+h_{3}(Y_{8},\bm{c}^{\{8\}})-3f_0.
\end{align}
Eq. (\ref{Eq: example-no-var}) thus shows that the $8$-dimensional problem has been reduced to a maximum dimensionality of $r=5$ using the first level check. Criterion (\ref{non-interaction}) is then tested on the $r$ (= 5) dimensional system with ${5}\choose{2}$ $=10$ cases. The set of interacting pairs of dimensions obtained from the interaction check is given by $I=\{\{1,2\},\{1,3\},\{2,3\},\{1,4\}\}$. Using information from the set $I$, $R_1=\{1,2,\dots,5\}$ is reduced to the following dimension set $Q$:
\begin{align}
Q&=\{Q_1, Q_2, Q_3\}=\{\{1,2,3\},\{1,4\},\{5\}\} \nonumber
\end{align}
We note that the presence of the 3-dimensional interaction $\{1,2,3\}$ have been derived from the interacting pairs $\{1,2\},\{1,3\}$ and $\{2,3\}$. This is how higher level interactions are derived from pairwise interaction results. Dimension set $S$ will then be given by
\begin{align}
S&=\{S_1, S_2, S_3, S_4, S_5, S_{6}\} \nonumber\\
&=\{\{6\},\{7\},\{8\},\{1,2,3\},\{1,4\},\{5\}\}\nonumber
\end{align}
Let $T$ be a collection of sets, which are the non-empty intersections between $S_i$ and $S_j$. We then have
\begin{equation}
T=\{\{1\}\} \nonumber\\
\end{equation} 
with $U=[1]$ and $V=4$. The function $g_{0}()$ will now be given by:
\begin{align}
g_{0}(Y_1,Y_2,....Y_{5}) &=h_{4}(Y_1,Y_2,Y_3,\bm{c}^{\{1,2,3\}}) +h_{5}(Y_1,Y_4,\bm{c}^{\{1,4\}}) \nonumber \\
& \ \ \ \ \ +h_{6}(Y_5,\bm{c}^{\{5\}}) 
 - p_{1}(Y_1,\bm{c}^{\{1\}})-f_0\nonumber\\
&=h_{4}(\bm{Y}_{S_{4}},\bm{c}^{S_{4}})
+h_{5}(\bm{Y}_{S_{5}},\bm{c}^{S_{5}}) \nonumber\\
& \ \ \ \ \ +h_{6}(\bm{Y}_{S_{6}},\bm{c}^{S_{6}}) 
 - p_{1}(\bm{Y}_{T_{1}},\bm{c}^{T_{1}})-f_0
\end{align}
Thus function $f(\bm{Y})$ is given by:
\begin{align}
f(\bm{Y})&=h_{4}(\bm{Y}_{S_{4}},\bm{c}^{S_{4}})
+h_{5}(\bm{Y}_{S_{5}},\bm{c}^{S_{5}}) +h_{6}(\bm{Y}_{S_{6}},\bm{c}^{S_{6}}) 
 - p_{1}(\bm{Y}_{T_{1}},\bm{c}^{T_{1}})-f_0\nonumber\\
& \ \ \ \ \ \ \ + h_{1}(Y_6,\bm{c}^{\{6\}})+h_{2}(Y_7,\bm{c}^{\{7\}})+h_{3}(Y_8,\bm{c}^{\{8\}}) -3f_0 \nonumber\\
&=h_{1}(\bm{Y}_{S_{1}},\bm{c}^{S_{1}})
+h_{2}(\bm{Y}_{S_{2}},\bm{c}^{S_{2}}) +h_{3}(\bm{Y}_{S_{3}},\bm{c}^{S_{3}}) 
+h_{4}(\bm{Y}_{S_{4}},\bm{c}^{S_{4}}) 
\nonumber\\
& \ \ \ \ \ \ \ +h_{5}(\bm{Y}_{S_{5}},\bm{c}^{S_{5}})+h_{6}(\bm{Y}_{S_{6}},\bm{c}^{S_{6}})  - p_{1}(\bm{Y}_{T_{1}},\bm{c}^{T_{1}}) -4f_0 \nonumber\\
& =\sum_{i=1}^{6}h_i(\bm{Y}_{S_i},\bm{c}^{S_i})
-\sum_{j=1}^{1}p_j(\bm{Y}_{T_j},\bm{c}^{T_j})-4f_0
\end{align}

\subsubsection{gPC Approximation Error}\label{cr: gpcerror}
Let us consider a $d$-dimensional domain where $1 \le d \le n$. Let $\bm{\xi_a}=\{\bm{\xi}_{\bm{a},1},\bm{\xi}_{\bm{a},2},\dots,\bm{\xi}_{\bm{a},m}\}$ be an array of $m$ Clenshaw-Curtis sparse grid points in dimension $d$ of depth level 2. There may also exist an additional array of $q$ unstructured points $\bm{\xi_b}=\{\bm{\xi}_{\bm{b},1},\bm{\xi}_{\bm{b},2},\dots,\bm{\xi}_{\bm{b},q}\}$ which have been previously evaluated. They correspond to sparse grid points in all ``predecessor" elements that are contained in the current domain. Let $\bm{u}_p$ be the second-order gPC approximation for the current domain corresponding to input points $\bm{\xi}$ where the gPC coefficients are calculated by solving a least squares problem given by Eq. (\ref{gPC_leastsq}) such that $\bm{\xi}=\{\bm{\xi_a},\bm{\xi_b}\}$ and $q+m=M$. Assuming $\bm{u}$ is the corresponding exact solution vector, the domain can be suitably approximated by the second-order gPC approximation if
\begin{align}\label{gPC-domain-approx}
\|\bm{u}^p-\bm{u}\|_{\infty}<\epsilon_1
\end{align}
If criterion (\ref{gPC-domain-approx}) is not satisfied, the domain is further subdivided into smaller elements along the center of its two most critical dimensions.
\subsection{Numerical Implementation}
The proposed algorithm is discussed below:\\
{\it{\textbf{Initialization and stopping criteria}}}. The dimension $n$ of the problem is first determined by the number of input random parameters considered in the model problem. $N_{iter}$ is the maximum number of iterations in the adaptive mesh refinement algorithm.
$V_{min}$ is a minimum hyper-volume fraction of the non-converged elements below which the subdivision into smaller elements is stopped.  When $N_{iter}$ is reached or the total hyper-volume fraction of the non-converged elements is less than $V_{min}$, the remaining non-converged elements are approximated by a first order gPC expansion and the algorithm terminates. Error tolerance parameters $\epsilon_1$ and $\epsilon_2$ are related to criteria (\ref{1D-check}), (\ref{no-variation}), (\ref{non-interaction}) and (\ref{gPC-domain-approx}). With decrease in the values of the chosen tolerance parameters, the approximation error also has a decreasing trend but with an increase in the computational cost because of more number of full model evaluations. \\
\\
{\it{\textbf{Checking global smoothness and possible dimensionality reduction}}}. This step initiates with the implementation of a first order gPC approximation in the original $n$-dimensional input space. The gPC coefficients are evaluated using the discrete projection method given by Eq. (\ref{discrete-proj}) using Clenshaw-Curtis sparse grid points of depth level $1$. The accuracy of the approximation is tested using criterion (\ref{gPC-domain-approx}). If the criterion is not satisfied, we go to the step of performing a one-dimensional (1-D) abrupt variation check. Otherwise, the first order gPC approximation is considered satisfactory and the algorithm skips to the surrogate value extraction step.\\
The 1-D abrupt variation check is now performed on the input domain to identify the influence of each dimension towards the output of interest. Criterion (\ref{1D-check}) is used to identify the critical dimensions while criterion (\ref{no-variation}) helps to reduce the $n$-dimensional problem to a number of problems with a maximum of $r$ dimensions where $r < n$. The interaction check is performed next, again on the global input domain using criterion
(\ref{non-interaction}) to further reduce the maximum dimensionality to $w ( <r)$ where $w=\max(|S_i|), \forall S_i \in S$.\\
If any of the dimensions are found to be critical based on the criterion of global abrupt variation, we directly go to the step of adaptive mesh refinement. Otherwise, a second order gPC approximation is now performed in the original $n$-dimensional input space using the discrete projection method. The function at the Clenshaw-Curtis sparse grid points of depth level 2 used for this approximation has already been evaluated in previous step of interaction check. Therefore, there is no extra computational cost involved for function evaluations in this step. The accuracy of the approximation is tested using criterion (\ref{gPC-domain-approx}). If the criterion is satisfied, the second order gPC approximation is considered satisfactory and the algorithm skips to the surrogate value generation step. Otherwise, we go to the next step.\\
\begin{algorithm}
\begin{algorithmic}
\caption{: \textbf{Summarized steps}}
\STATE \textbf{\underline{Initialization}}\\
\STATE Set $n$, $N_{iter}$, $V_{min}$, $\epsilon_1$, and $\epsilon_2$.\\
\textbf{\underline{Global checks and dimensionality reduction}}\\

\STATE perform first order gPC approximation using Eq. (\ref{discrete-proj})

\IF { $||\bm{u}^p-\bm{u}||_{\infty}<\epsilon_1$ (see Eq. (\ref{gPC-domain-approx}))}

\STATE go to the \underline{Surrogate value extraction} step

\ELSE

\STATE perform abrupt variation check using criterion (\ref{1D-check})

\STATE perform dimensionality reduction using criteria (\ref{no-variation}) and (\ref{non-interaction}) to form lower dimensional sub-problems.

\IF {$||\bm{u}_p^{(i)}-\bm{u}^{(i)}||_{\infty}<\epsilon_1$ (see Eq. (\ref{1D-check})) for all dimensions}

\STATE perform second order gPC approximation using Eq. (\ref{discrete-proj})

\IF { $||\bm{u}^p-\bm{u}||_{\infty}<\epsilon_1$ (see Eq. (\ref{gPC-domain-approx}))}

\STATE go to the \underline{Surrogate value extraction} step

\ELSE

\STATE go to the \underline{Adaptive mesh refinement} step

\ENDIF

\ELSE

\STATE go to the \underline{Adaptive mesh refinement} step

\ENDIF\\

\ENDIF\\
\end{algorithmic}
\end{algorithm}

\begin{algorithm}
\begin{algorithmic}
\caption*{\textbf{Algorithm 1 : Summarized steps (continued)}}
\STATE
\textbf{\underline{Adaptive mesh refinement}}\\
\FORALL {sub-dimensional problems}
	\STATE check abrupt variations using criterion (\ref{1D-check})
	
	\IF {criterion (\ref{1D-check}) is satisfied}
		\STATE check gPC approximation using criterion (\ref{gPC-domain-approx})
		
		\IF {criterion (\ref{gPC-domain-approx}) is not satisfied} 
			\STATE subdivide the element along the center of its two most critical dimensions
 		\ENDIF
		
	\ELSE
		
		\STATE subdivide the element along the center of its two most critical dimensions
		 	\ENDIF
 \ENDFOR\\
 \textbf{\underline{Surrogate value extraction}}
 \STATE extract output values corresponding to query inputs from the approximate model obtained.
\end{algorithmic}
\end{algorithm}
{\it{\textbf{Adaptive mesh refinement}}}. This part of the algorithm in general deals with $(N_S+N_T)$ low dimensional subproblems as mentioned in section \ref{cr: interaction}. For a subproblem $\euscr{P}_i \ ( 1 \le i \le N_S+N_T) \} )$, the algorithm initiates with the subdivision of the original domains into elements along its two most critical dimensions. The iteration count $Iter$ starts here. For each of the $E_{\euscr{P}_i}$ elements formed in $\euscr{P}_i$ in a certain iteration, an abrupt variation check is performed as was done on the original $n$-dimensional domain. If the second-order approximation criterion (\ref{1D-check}) is not met, the element $E_{\euscr{P}_i}^j \ (j \in \{1,2,...,E_{\euscr{P}_i} \})$  is again subdivided into subelements along its two most critical dimensions. Satisfaction of criterion (\ref{1D-check}) implies there are no 
abrupt variations in the current element. This leads to checking criterion (\ref{gPC-domain-approx}) for second order gPC approximation in the whole element. If that criterion is met, the element $E_{\euscr{P}_i}^j$ is said to have converged for the given tolerance $\epsilon_1$ and can be suitably approximated by a second order gPC approximation. The polynomial order, the coefficient vector and the range of the converged element is then stored for future surrogate retrieval. If criterion (\ref{gPC-domain-approx}) is not satisfied, then the element is also subdivided into smaller elements. This procedure is performed for all $E_{\euscr{P}_i}$ elements and all the new subelements formed undergo similar operations at the next iteration $Iter=Iter+1$. At the end of each iteration, the hyper volume $V$ of the subelements created and the number of iterations $Iter$ are compared with the corresponding critical values $V_{min}$ and $N_{iter}$ respectively to check if either of the two stopping criteria is met. If the stopping condition gets satisified, then all the remaining subelements are approximated by a first order gPC approximation. After meeting the stopping criteria, the next subproblem is taken up and we repeat the process of characterizing it.\\

{\it{\textbf{Surrogate value extraction}}}. After having characterized the $n$-dimensional problem through the various steps mentioned, the final step is to generate output values corresponding to arbitrary query input points in the $n$-dimensional domain and also output statistics, such as, mean. Output value estimation corresponding to a query input involves locating the element in which the query point lies in each subproblem. The stored information for that element is then retrieved to generate the local surrogate output values in each subproblem, which are then combined together to get the global output value. Mean value estimation is performed by evaluating the integration in each of the elements in each subproblem. For each subproblem, the global mean is calculated by the weighted average of local means corresponding to each element, and the weight is based on the ratio of the hyper-volume of the elements and the hyper-volume of the whole domain.\\
A summary of the all the above steps is given in Algorithm 1.

\section{Numerical Results}
In this section, SCAMR is applied to a variety of functions with smoothness as well as discontinuities and input dimensions as high as 300. Through these examples, its performance is tested against existing efficient algorithms, like, ASGC \cite{ma2009adaptive}, HDMR-ASGC \cite{ma2010adaptive} and MEPCM-A \cite{foo2010multi}.
\subsection{Demonstration of SCAMR Performance}
We first demonstrate the effectiveness and efficiency of the proposed SCAMR method using simple smooth functions with random input spaces of different dimensions. Then, we will focus on functions with non-smoothness or discontinuities in random space, as well as a high-dimensional stochastic elliptic problem. Our results are compared to those from ASGC method since both approaches use low order polynomials as a basis and both use adaptivity to track discontinuities. Specifically, we compare the root mean squared error calculated using 
$N=10^5$ randomly generated samples, given by 
\begin{equation}
\epsilon=\sqrt{\frac{1}{N}\sum_{i=1}^N (f(\mathbf{x_i})-\tilde{f}(\mathbf{x_i}))^2},
 \end{equation}
where $f$ is the exact function and $\tilde{f}$ is the numerical approximation using ASGC or SCAMR. 

\subsubsection{Performance of SCAMR on Smooth Functions}
\begin{figure}
\centering
\begin{subfigure}[b]{0.4\textwidth}
\centering
\includegraphics[width=1.2\linewidth, height=4cm]{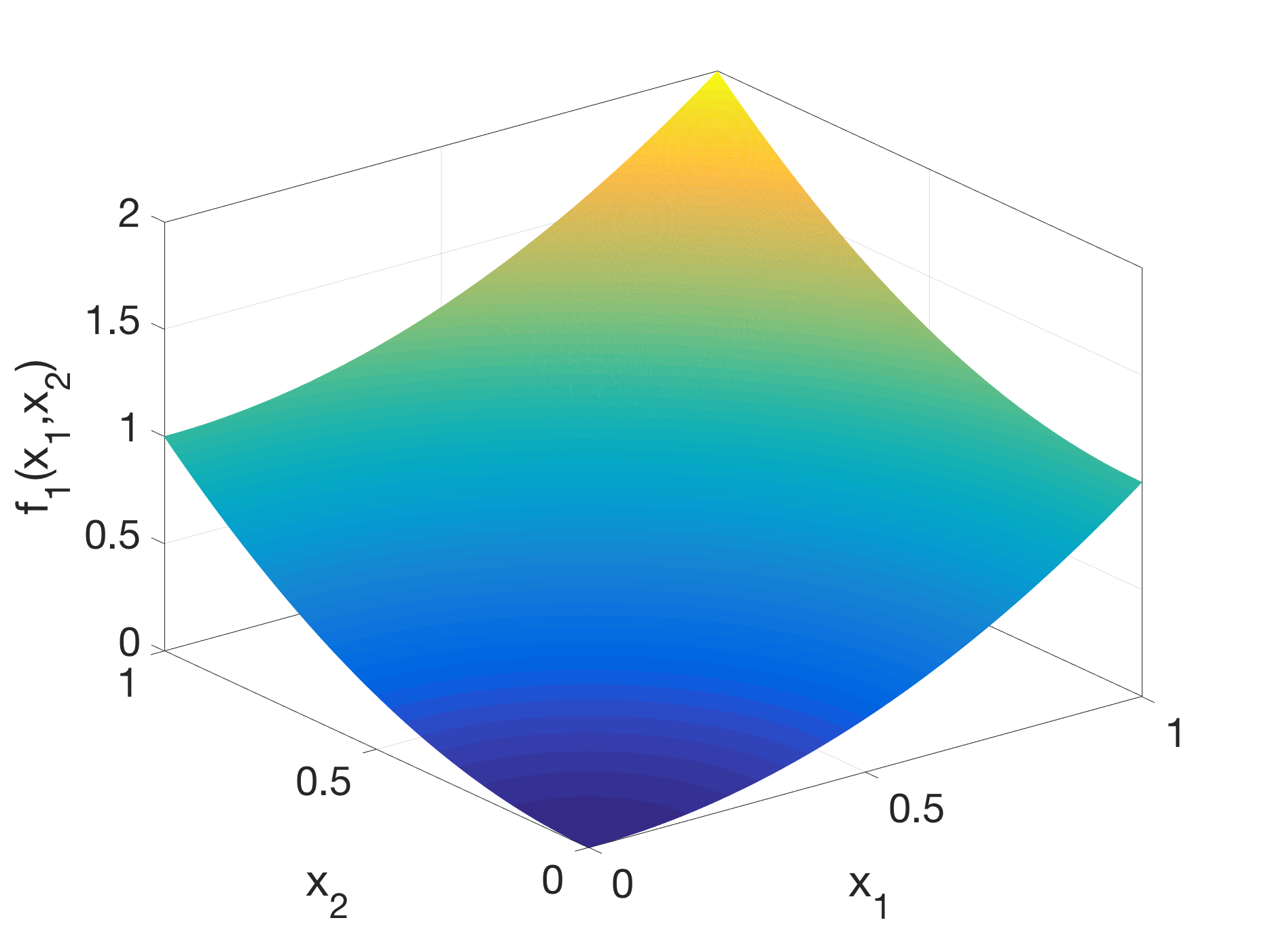} 
\caption{} 
\label{fig:subim131}
\end{subfigure}
\hfill
\begin{subfigure}[b]{0.4\textwidth}
\centering
\includegraphics[width=1.2\linewidth, height=4cm]{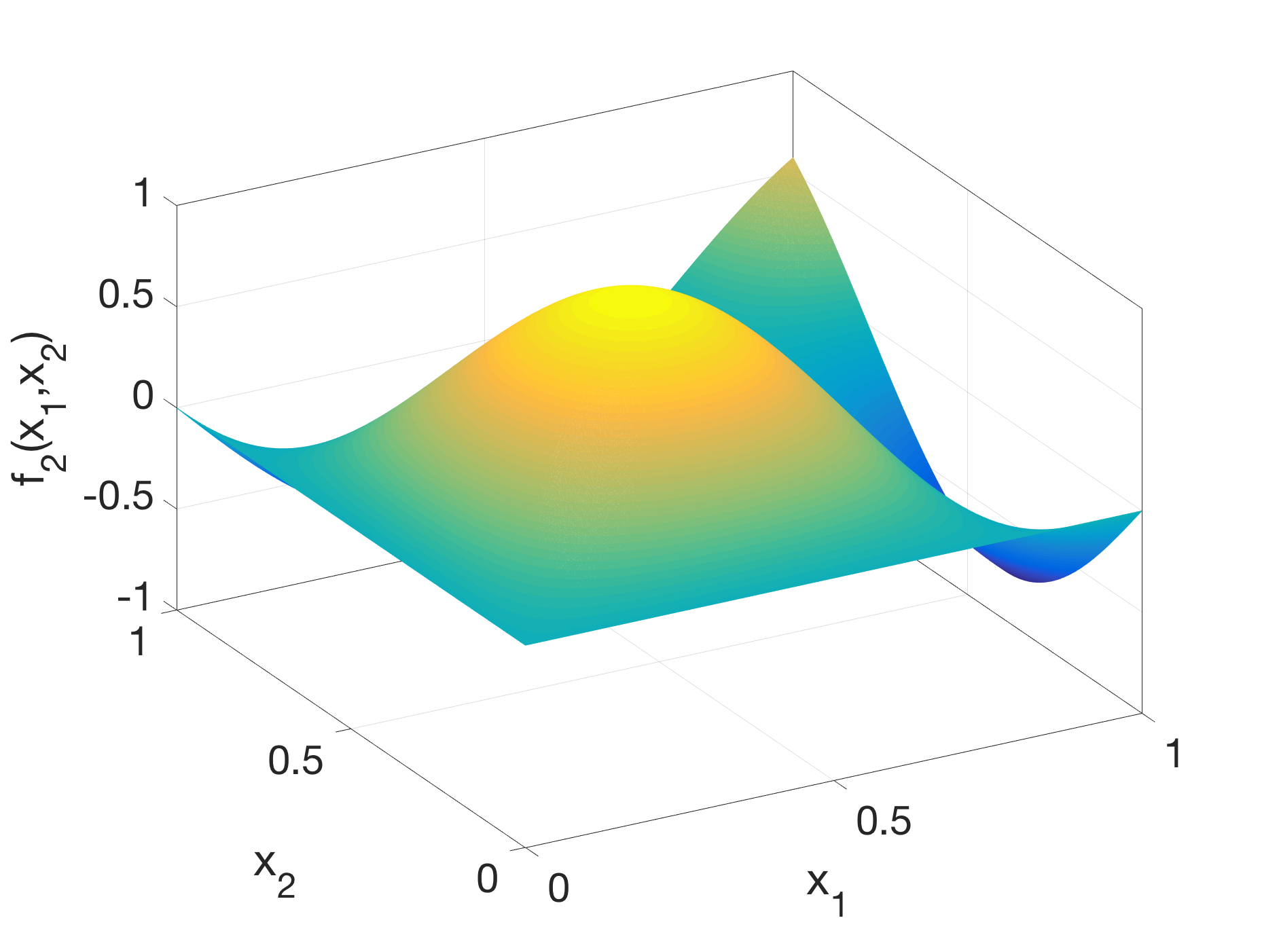} 
\caption{} 
\label{fig:subim131}
\end{subfigure}
\hfill
\begin{subfigure}[b]{0.4\textwidth}
\includegraphics[width=1.2\linewidth, height=4.4cm]{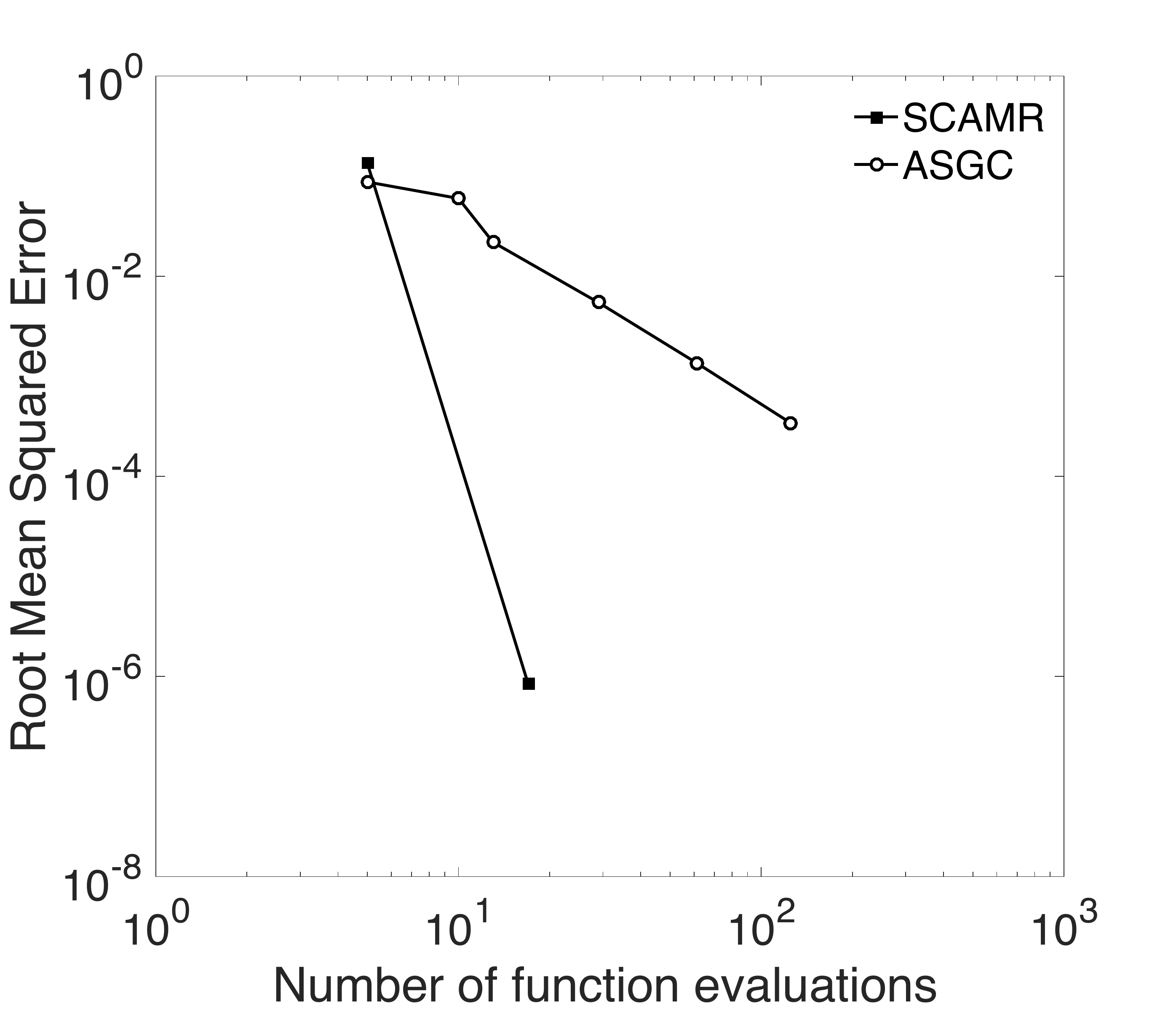} 
\caption{} 
\label{fig:subim131}
\end{subfigure}
\hfill
\begin{subfigure}[b]{0.4\textwidth}
\includegraphics[width=1.2\linewidth, height=4.4cm]{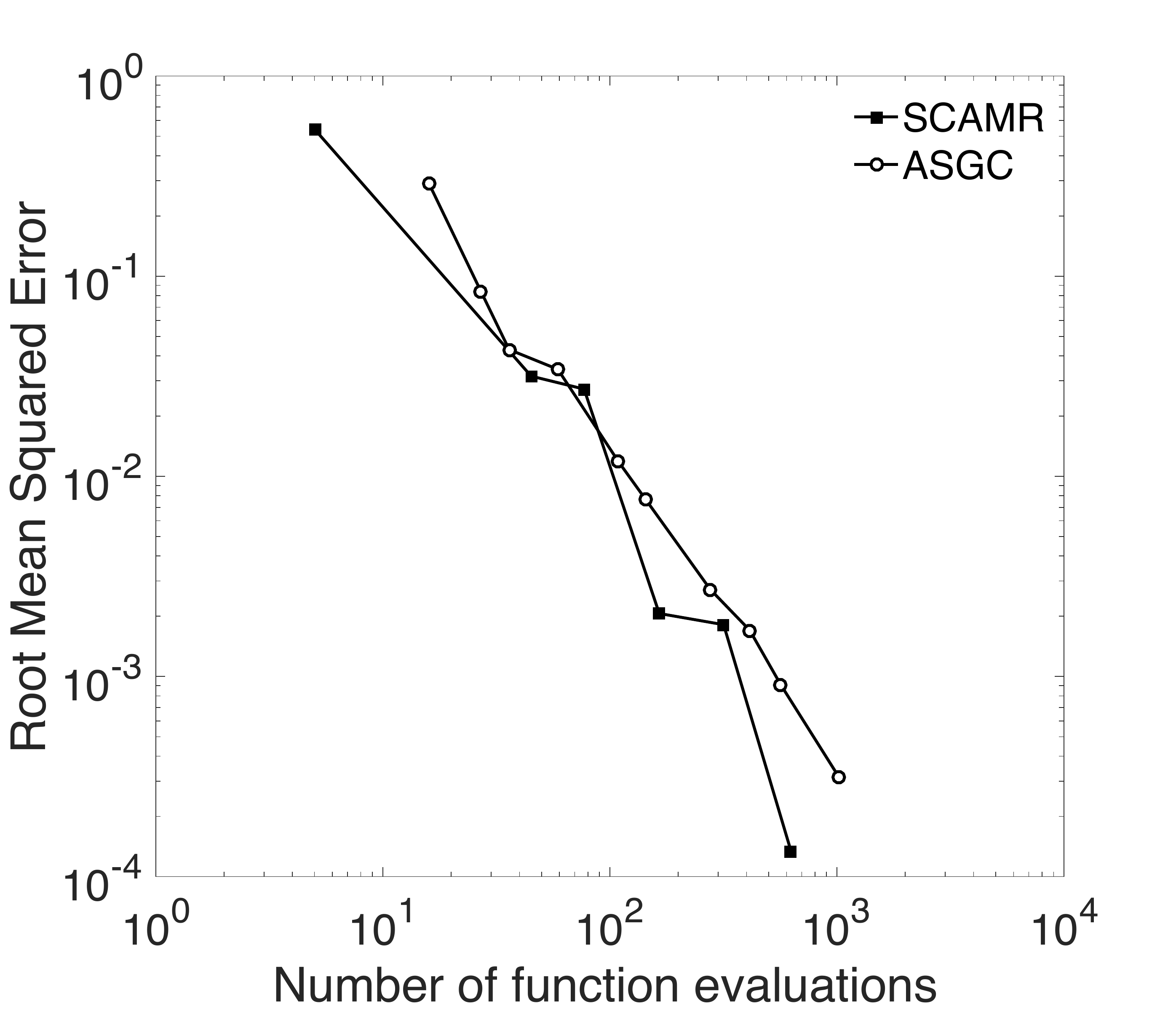} 
\caption{} 
\label{fig:subim131}
\end{subfigure}
\caption{Results for 2D smooth functions: (a) exact function output for $f_1$; (b) exact function output for $f_2$; (c) error of estimated $f_1$ using SCAMR and the ASGC method; and (d) error of estimated $f_2$ using SCAMR and the ASGC method.}
\label{fig:image3}
\end{figure}

We first implement the proposed method on a few simple smooth functions with random inputs in different dimensions.
The two-dimensional test functions are quadratic and sine functions defined as follows.
\begin{eqnarray}
f_1(x_1,x_2)&=&x_1^2+x_2^2, \\
f_2(x_1,x_2)&=&\sin(4x_1)\sin(4x_2),
\end{eqnarray}
where $x_i$ are i.i.d. uniform random variables in $[0,1]$ ($i=1,2$). The exact functions are provided in Fig. \ref{fig:image3}(a,b) for $f_1$ and $f_2$ respectively. Clearly, the product of sine functions $f_2$ exhibits more abrupt variations than the summation of quadratic functions $f_1$ in the $[0,1]^2$ domain; therefore, one would expect slower convergence of the numerical approximation for $f_2$. The numerical errors of SCAMR method are provided in Fig. \ref{fig:image3}(c,d), and compared to those from ASGC method. From the results, one can observe that i) both SCAMR and ASGC methods have slower convergence for $f_2$ compared to $f_1$ as we expected, and ii) our proposed SCAMR approach converges faster than ASGC for both the functions.

\begin{figure}
\centering
\begin{subfigure}[b]{0.4\textwidth}
\centering
\includegraphics[width=1.2\linewidth, height=4.4cm]{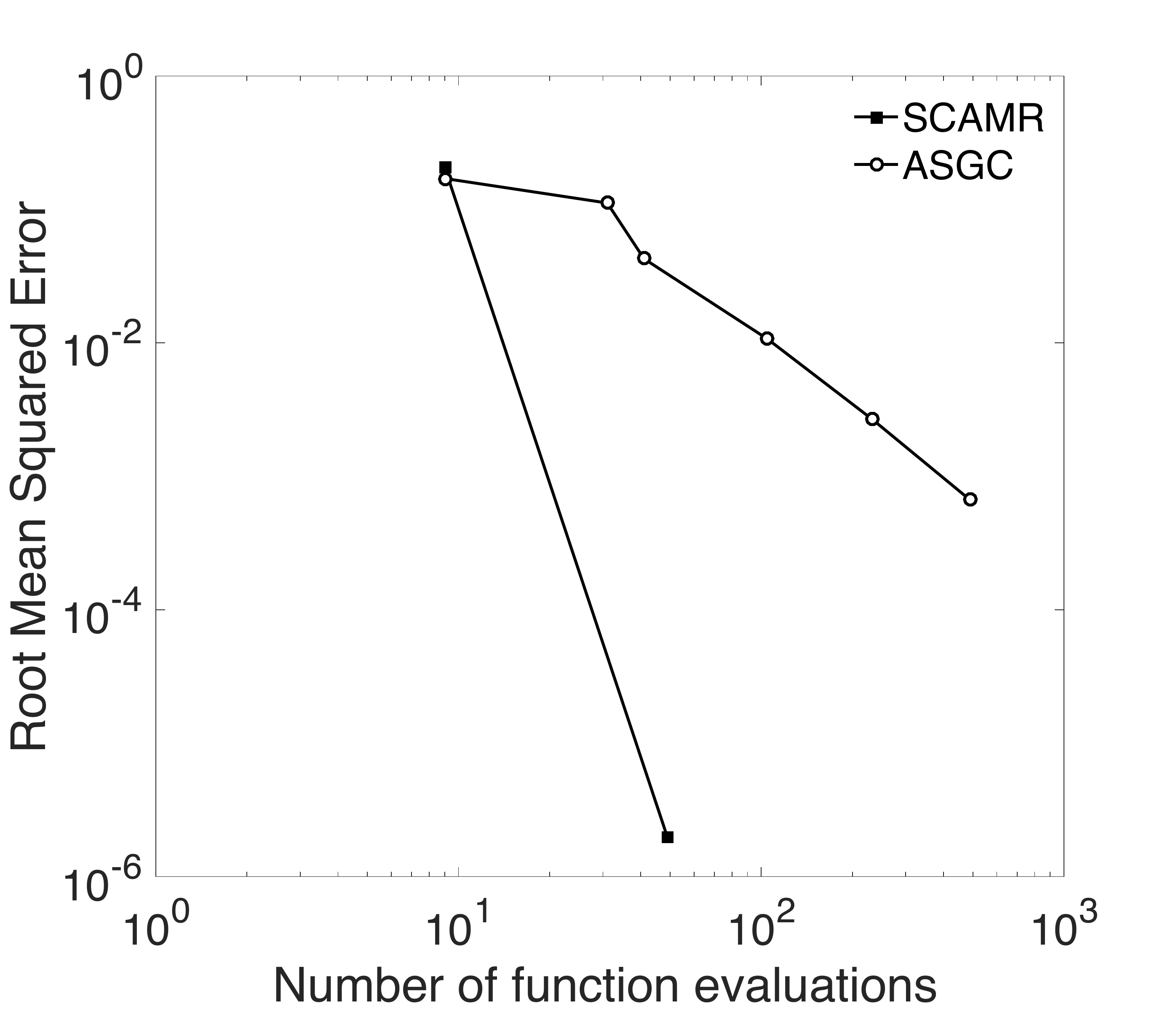} 
\caption{} 
\label{fig:subim131}
\end{subfigure}
\hfill
\begin{subfigure}[b]{0.4\textwidth}
\centering
\includegraphics[width=1.2\linewidth, height=4.4cm]{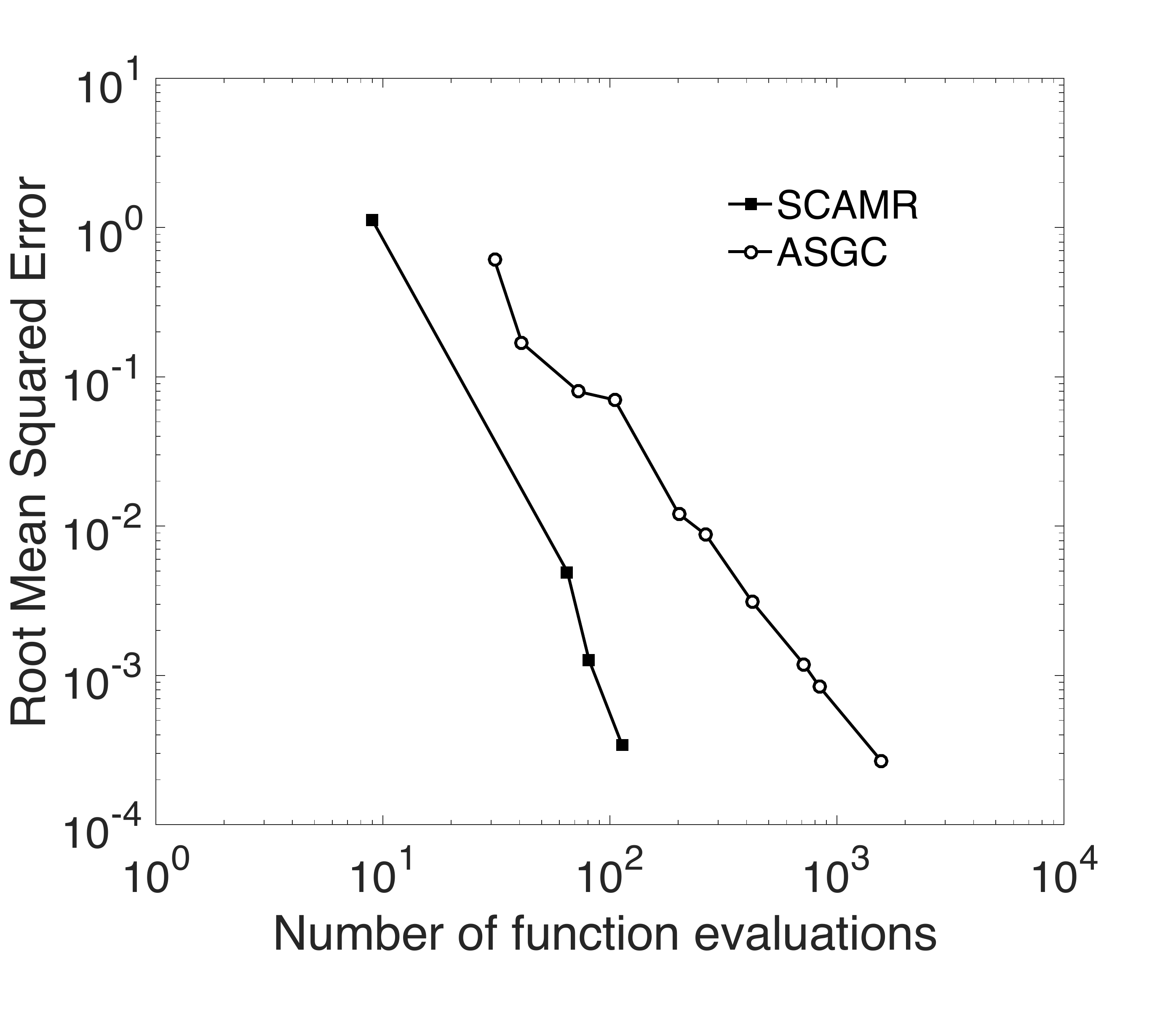} 
\caption{} 
\label{fig:subim131}
\end{subfigure}
\hfill
\begin{subfigure}[b]{0.4\textwidth}
\includegraphics[width=1.2\linewidth, height=4.4cm]{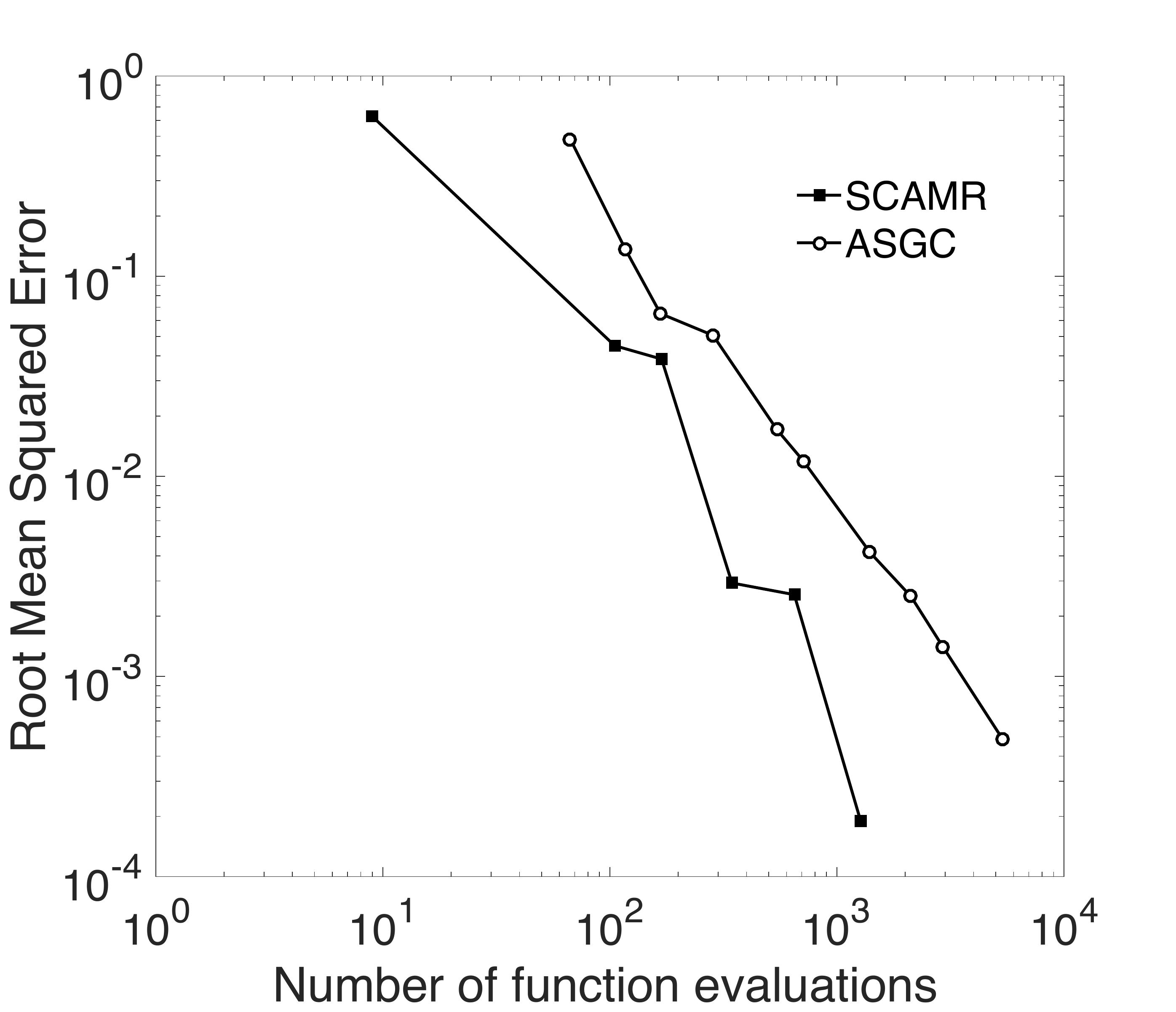} 
\caption{} 
\label{fig:subim131}
\end{subfigure}
\hfill
\begin{subfigure}[b]{0.4\textwidth}
\includegraphics[width=1.2\linewidth, height=4.4cm]{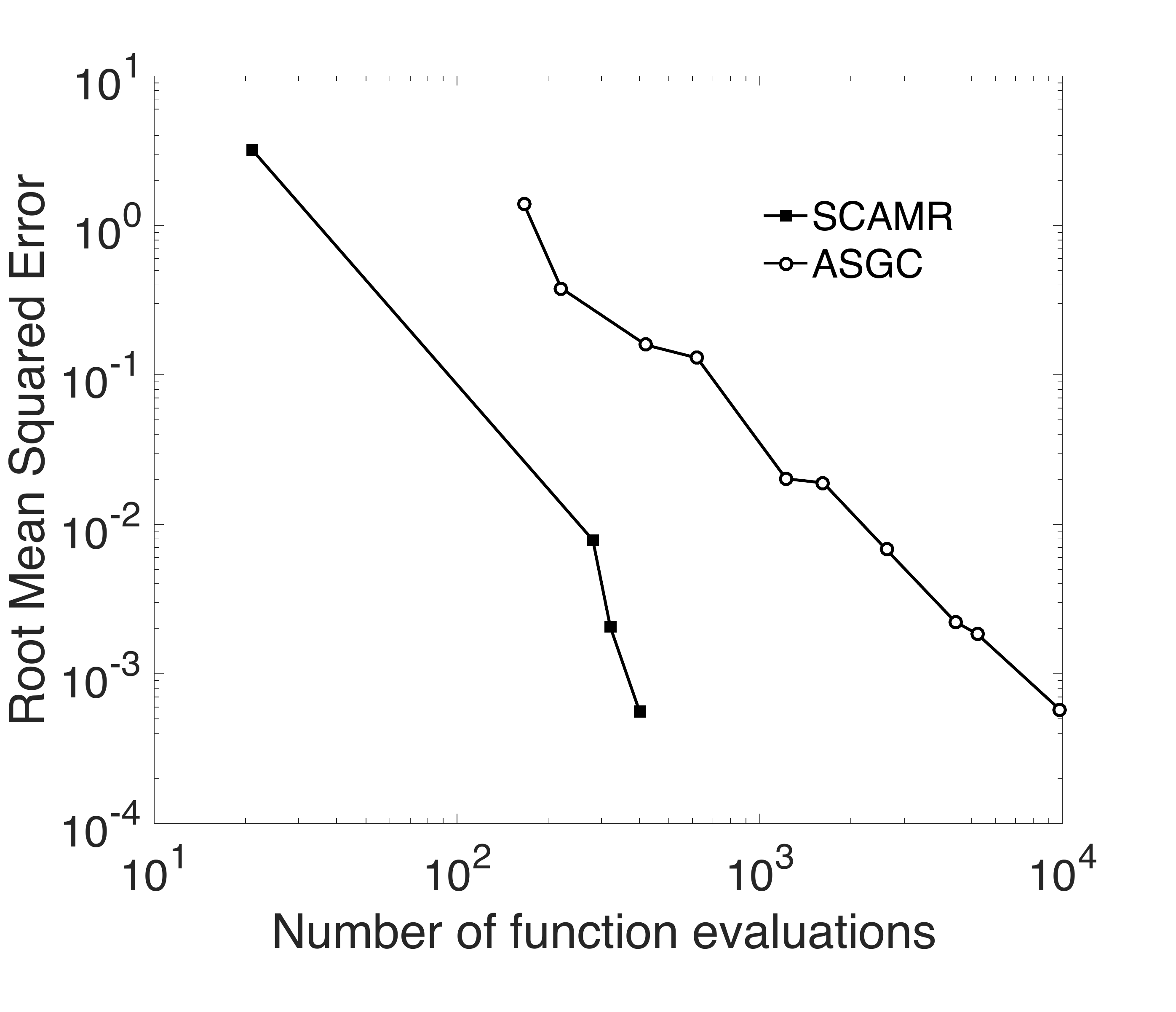} 
\caption{} 
\label{fig:subim131}
\end{subfigure}
\caption{Error analysis of SCAMR and ASGC methods for 4D and 10D smooth functions: (a) 4D $f_3$, (b) 4D $f_4$, (c) 4D $f_5$, and (d) 10D $f_6$.}
\label{fig:image4}
\end{figure}

We extend two-dimensional quadratic and sine functions to four and ten dimensions as follows.
\begin{eqnarray}
f_3(x_1,x_2,x_3,x_4)&=&\sum_{i=1}^4x_i^2, \\
f_4(x_1,x_2,x_3,x_4)&=&\sum_{i=1}^4\sin(4x_i),\\
f_5(x_1,x_2,x_3,x_4)&=&\sin(4x_1)\sin(4x_2)+\sin(4x_3)\sin(4x_4),\\
f_6(x_1,x_2, \hdots, x_{10})&=&\sum_{i=1}^{10}\sin(4x_i),  
\end{eqnarray}
where $x_i$ are i.i.d. uniform random variables in $[0,1]$ ($i=1,2,\hdots,10$). The functions $f_3$, $f_4$ and $f_6$ are independent of the interaction terms between the inputs, while $f_5$ depends on some interaction terms between the inputs. The numerical errors of both the SCAMR and the ASGC methods are provided in Fig. \ref{fig:image4} with respect to number of function evaluations. The numerical approximation from both methods converges slower as the complexity of the function increases, such as, from a polynomial function to a sine function, from an additive function to a multiplicative function or from a lower dimensional (4-D) function to a higher dimensional (10-D) function. Fig. \ref{fig:image4} shows that SCAMR converges faster than ASGC for all four smooth functions. 

Having tested the SCAMR approach on smooth functions with random inputs in different dimensions, we will next discuss its performance on non-smooth functions. 

\subsubsection{Performance of SCAMR on Functions with Line Singularity}
\begin{figure}
\centering
\includegraphics[width=0.6\linewidth, height=5cm]{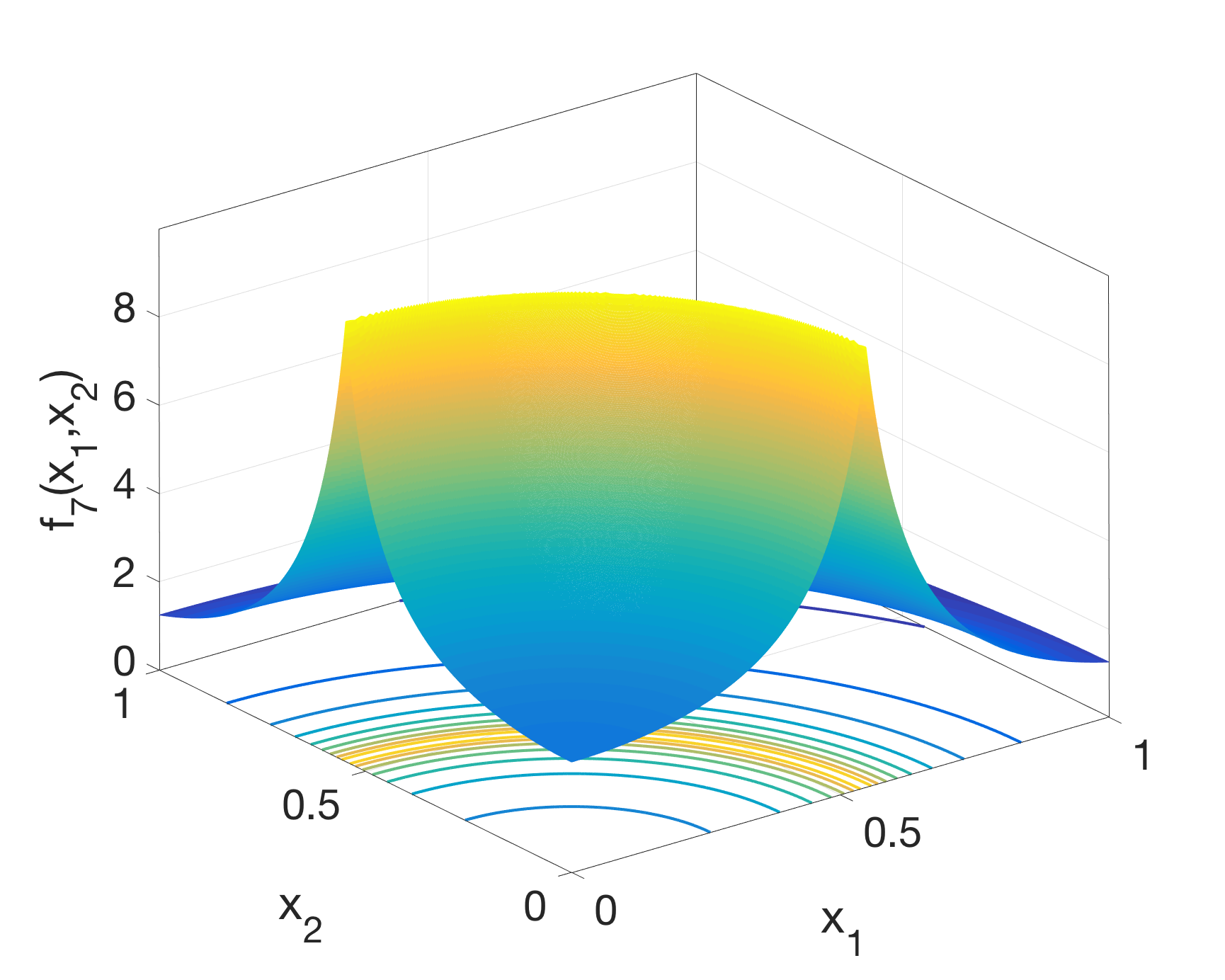} 
\caption{Surface plot of function $f_7(x_1,x_2)$.} 
\label{fig:image5}
\end{figure}

Here we adopt the same 2D function with line singularity as in \cite{ma2009adaptive}. 
\begin{equation}
f_7(x_1,x_2)=\frac{1}{|0.3-x_1^2-x_2^2|+0.1}.
\end{equation}
The function is plotted in Fig. \ref{fig:image5}. Clearly, the function has a $C^1$ discontinuity going across both $x_1$ and $x_2$ directions.

The 4D and 10D extensions of the above function are defined as
\begin{eqnarray}
f_8(x_1,x_2,x_3,x_4)&=&\frac{1}{|0.3-x_1^2-x_2^2|+0.1}+ \sum_{i=3}^4x_i,\\
f_9(x_1,x_2,\hdots, x_{10})&=&\frac{1}{|0.3-x_1^2-x_2^2|+0.1}+ \sum_{i=3}^{10}x_i
\end{eqnarray}
where $x_i$ are i.i.d. uniform random variables in $[0,1]$ ($i=1,2,\hdots,10$). Notice that the added dimensions in $f_8$ and $f_9$ are not interactive with $x_1$ and $x_2$. Therefore one would expect that the computational cost will not increase dramatically as the dimension increases.

\begin{figure}
\centering
\begin{subfigure}[b]{0.4\textwidth}
\centering
\includegraphics[width=1.2\linewidth, height=4cm]{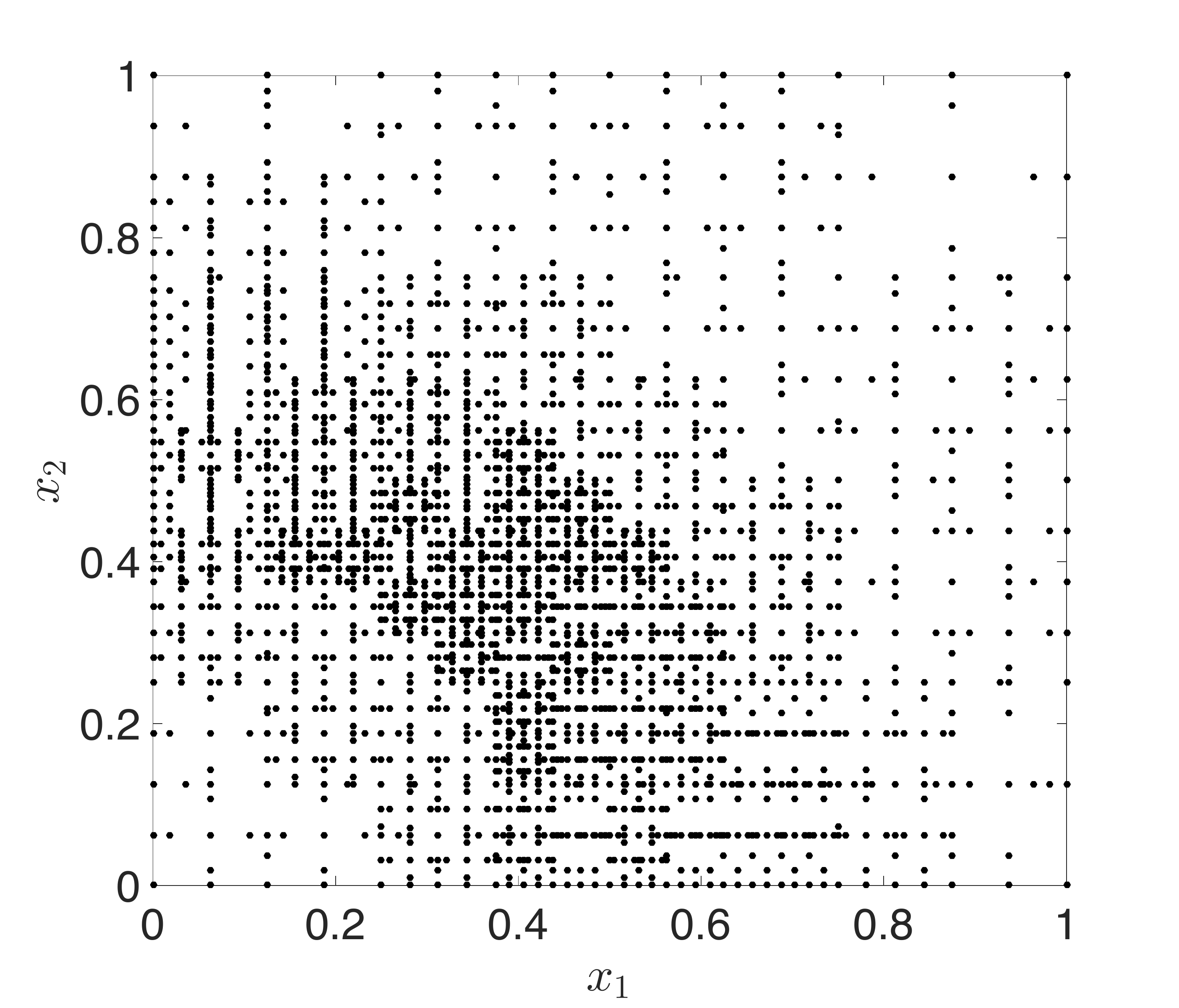} 
\caption{} 
\label{fig:subim131}
\end{subfigure}
\hfill
\begin{subfigure}[b]{0.4\textwidth}
\centering
\includegraphics[width=1.2\linewidth, height=4cm]{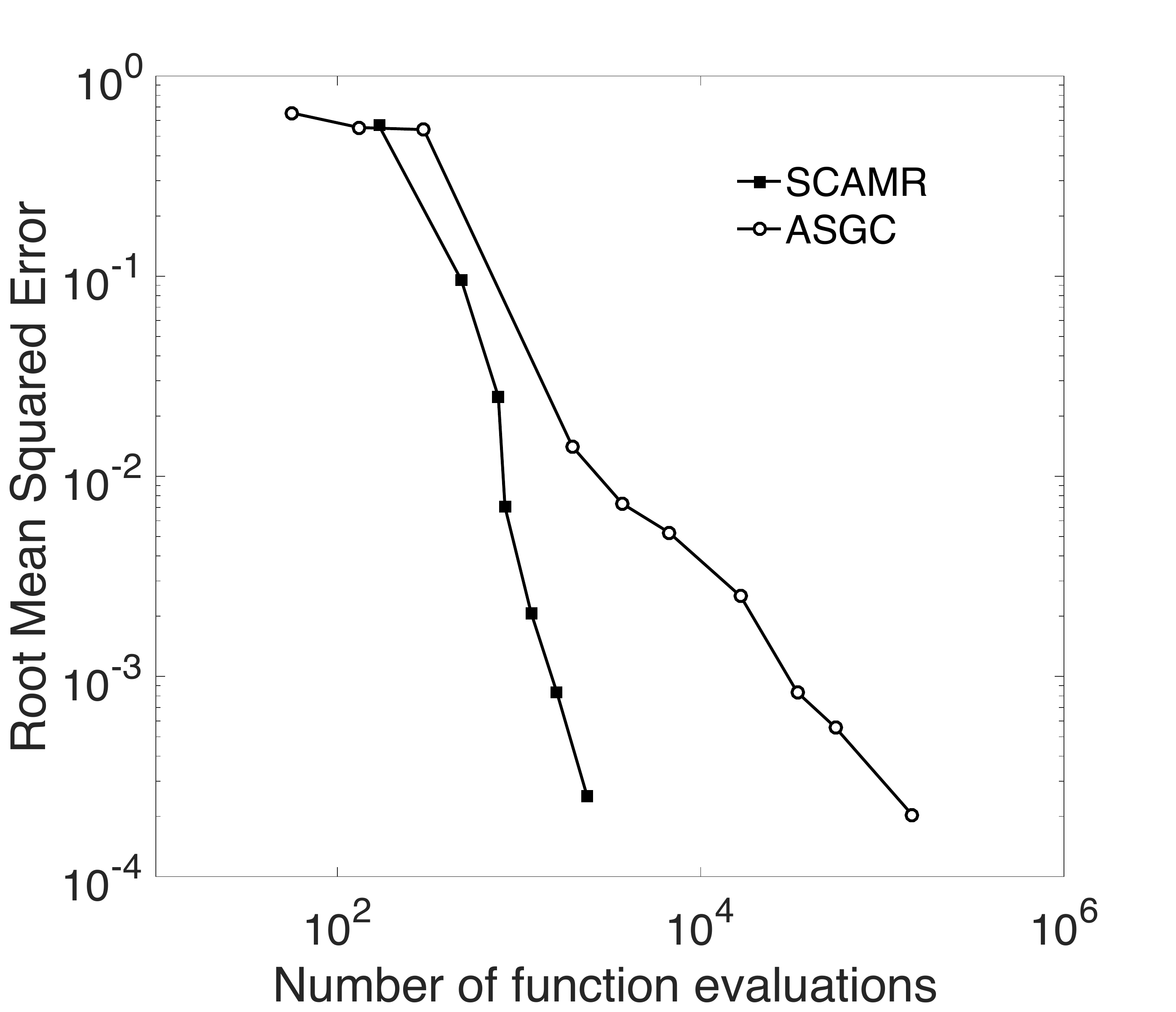} 
\caption{} 
\label{fig:subim131}
\end{subfigure}
\hfill
\begin{subfigure}[b]{0.4\textwidth}
\centering
\includegraphics[width=1.2\linewidth, height=4cm]{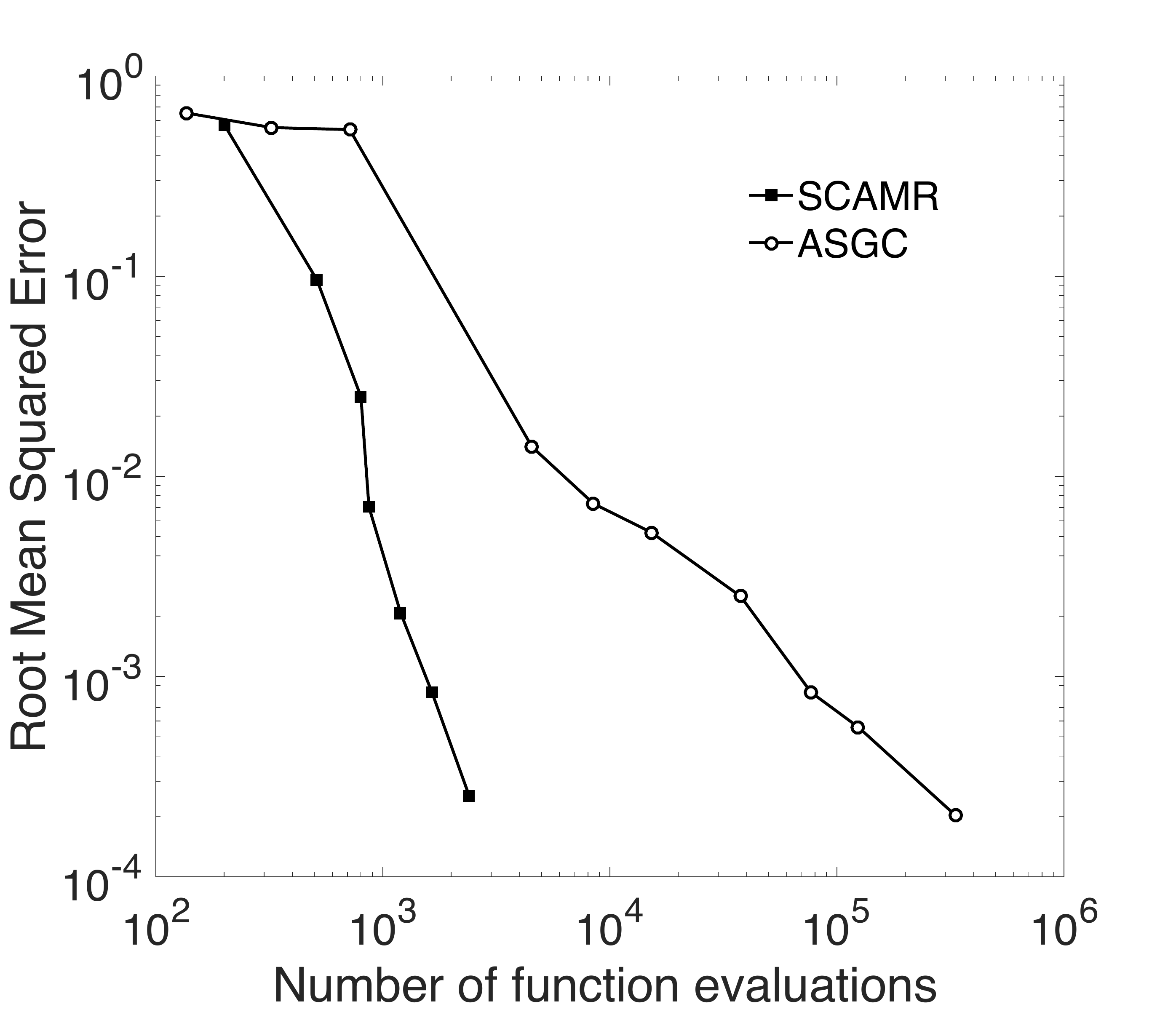} 
\caption{} 
\label{fig:subim131}
\end{subfigure}
\hfill
\begin{subfigure}[b]{0.4\textwidth}
\includegraphics[width=1.2\linewidth, height=4cm]{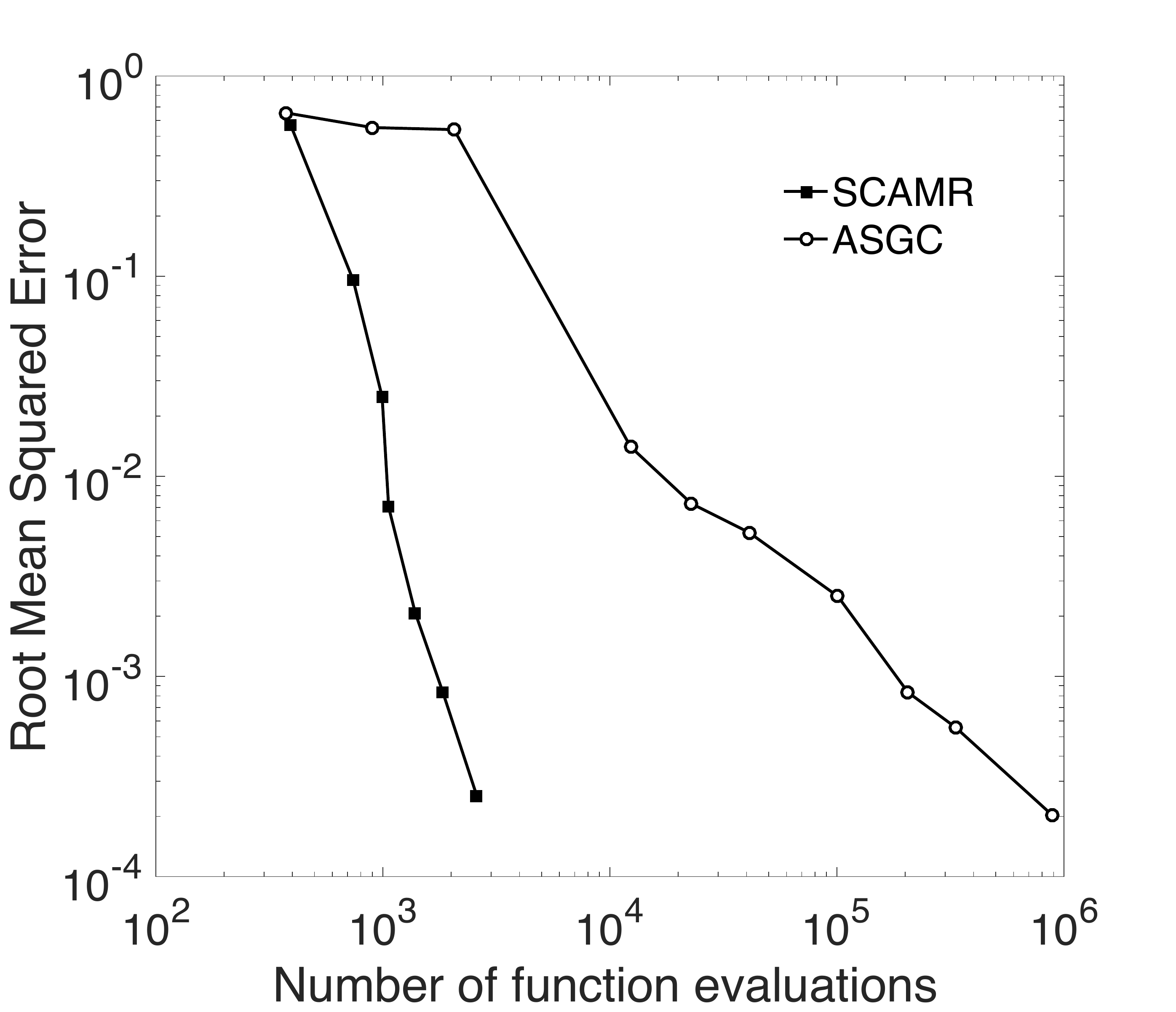} 
\caption{} 
\label{fig:subim131}
\end{subfigure}
\caption{Input domain and error analysis for functions with line singularity: (a) input domain for function $f_7$, (b) numerical error as a function of the number of samples for 2D $f_7$, (c) numerical error as a function of the number of samples for 4D $f_8$, and (d) numerical error as a function of the number of samples for 10D $f_9$. }
\label{fig:image6}
\end{figure}

The proposed SCAMR approach is implemented for the above 2-D, 4-D and 10-D functions. The locations of function evaluations for the 2-D function $f_7$ are plotted in Fig. \ref{fig:image6}a. The plot shows that the line singularity is well captured by the approach and more function evaluations are required in the area of line singularity as expected. The error analysis of the numerical approximations are provided in Fig. \ref{fig:image6}(b-d) for functions $f_7$, $f_8$ and $f_9$, respectively. From the figure, one can observe that the convergence rates of SCAMR are similar for the three functions with different dimensions as expected. The SCAMR approach converges faster than ASGC for all three functions. 

\subsubsection{Performance of SCAMR on Functions with $C^0$ discontinuity}
\begin{figure}
\centering
\includegraphics[width=0.6\linewidth, height=4cm]{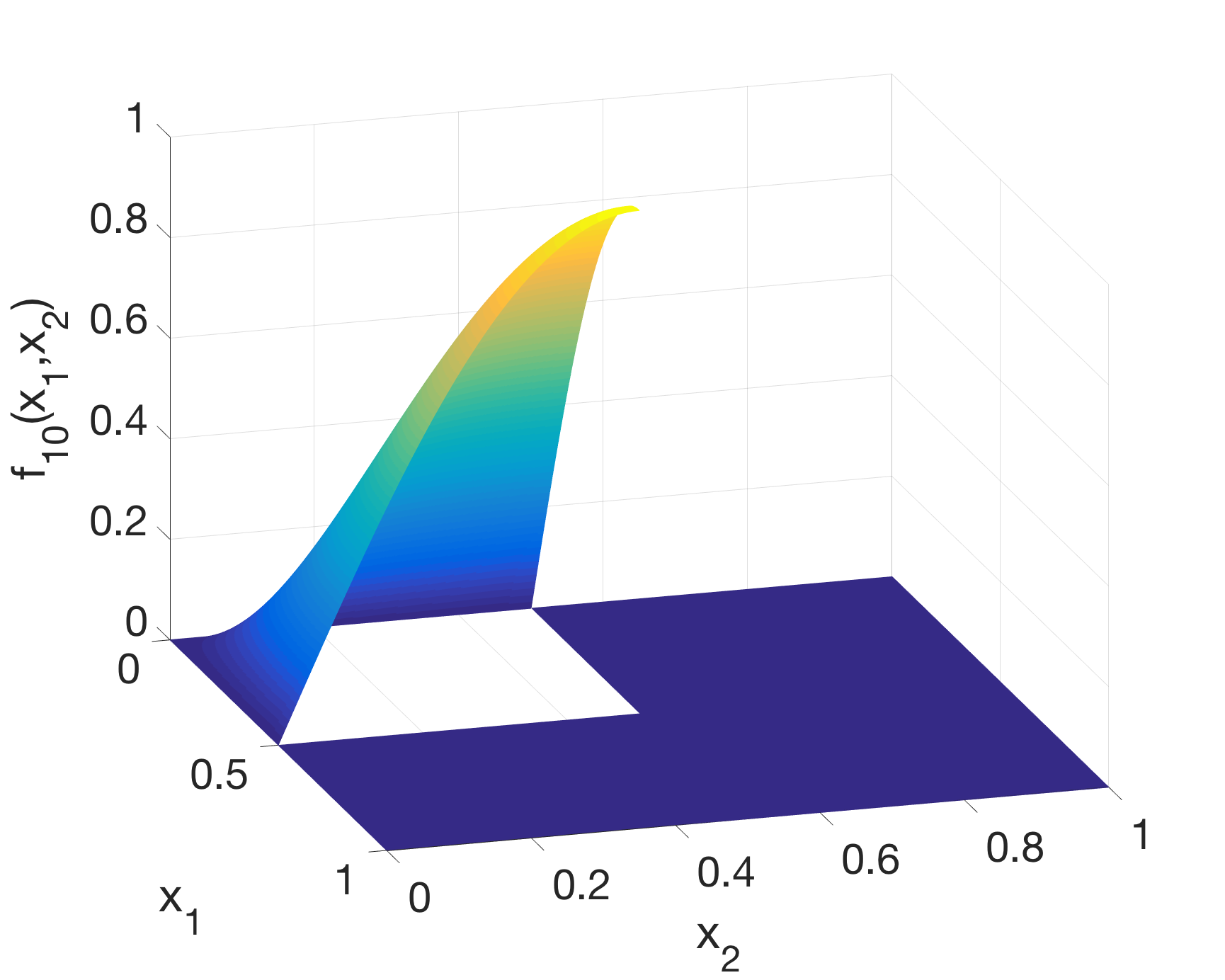} 
\caption{Surface plot of function $f_{10}(x_1,x_2)$.} 
\label{fig:image7}
\end{figure}

SCAMR is tested on another 2-D function, this one with a $C^0$ discontinuity as in  \cite{agarwal2009domain}:
\[
    f_{10}(x_1,x_2)= 
\begin{cases}
    0,& \text{if } x_1 \geq 0.5 \text{ or } x_2 \ge 0.5,\\
    \sin(\pi x_1)\sin(\pi x_2),              & \text{otherwise}
\end{cases}
\]
The function is plotted in Fig. \ref{fig:image7}.

Similarly, we extend it to 4-D and 10-D functions with discontinuity as
\[
    f_{11}(x_1,x_2,x_3,x_4)= 
\begin{cases}
    \sum_{i=3}^4x_i,& \text{if } x_1 \geq 0.5 \text{ or } x_2 \ge 0.5,\\
    \sin(\pi x_1)\sin(\pi x_2)+\sum_{i=3}^4x_i,              & \text{otherwise}
\end{cases}
\]
and
\[
    f_{12}(\textbf{x})= 
\begin{cases}
    \sum_{i=3}^{10}x_i,& \text{if } x_1 \geq 0.5 \text{ or } x_2 \ge 0.5,\\
    \sin(\pi x_1)\sin(\pi x_2)+\sum_{i=3}^{10}x_i,              & \text{otherwise}
\end{cases}
\]
where  $\textbf{x}=\{x_1,x_2,\ldots,x_{10}\}$.

\begin{figure}
\centering
\begin{subfigure}[b]{0.4\textwidth}
\centering
\includegraphics[width=1.2\linewidth, height=4cm]{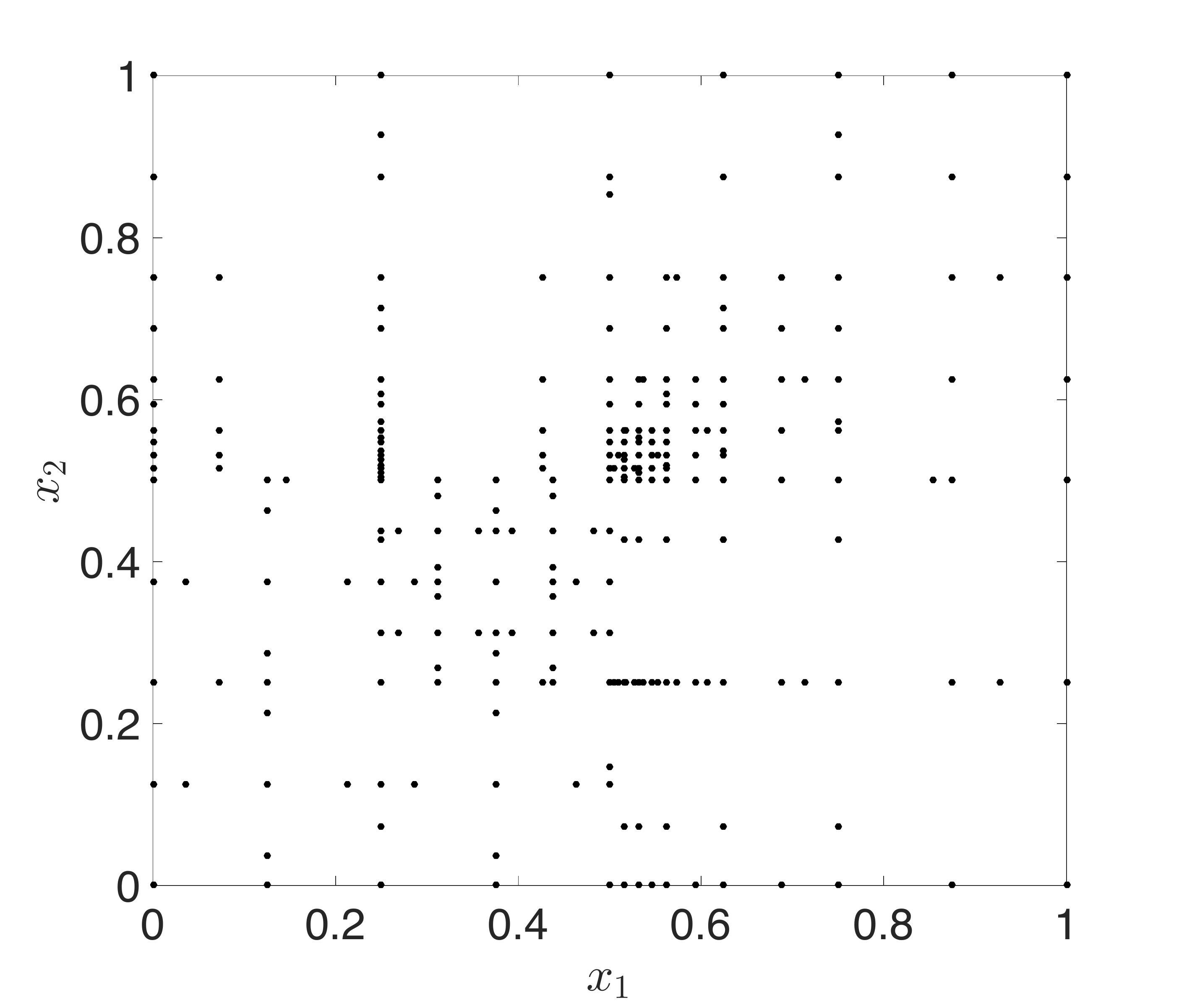} 
\caption{} 
\label{fig:subim131}
\end{subfigure}
\hfill
\begin{subfigure}[b]{0.4\textwidth}
\centering
\includegraphics[width=1.2\linewidth, height=4cm]{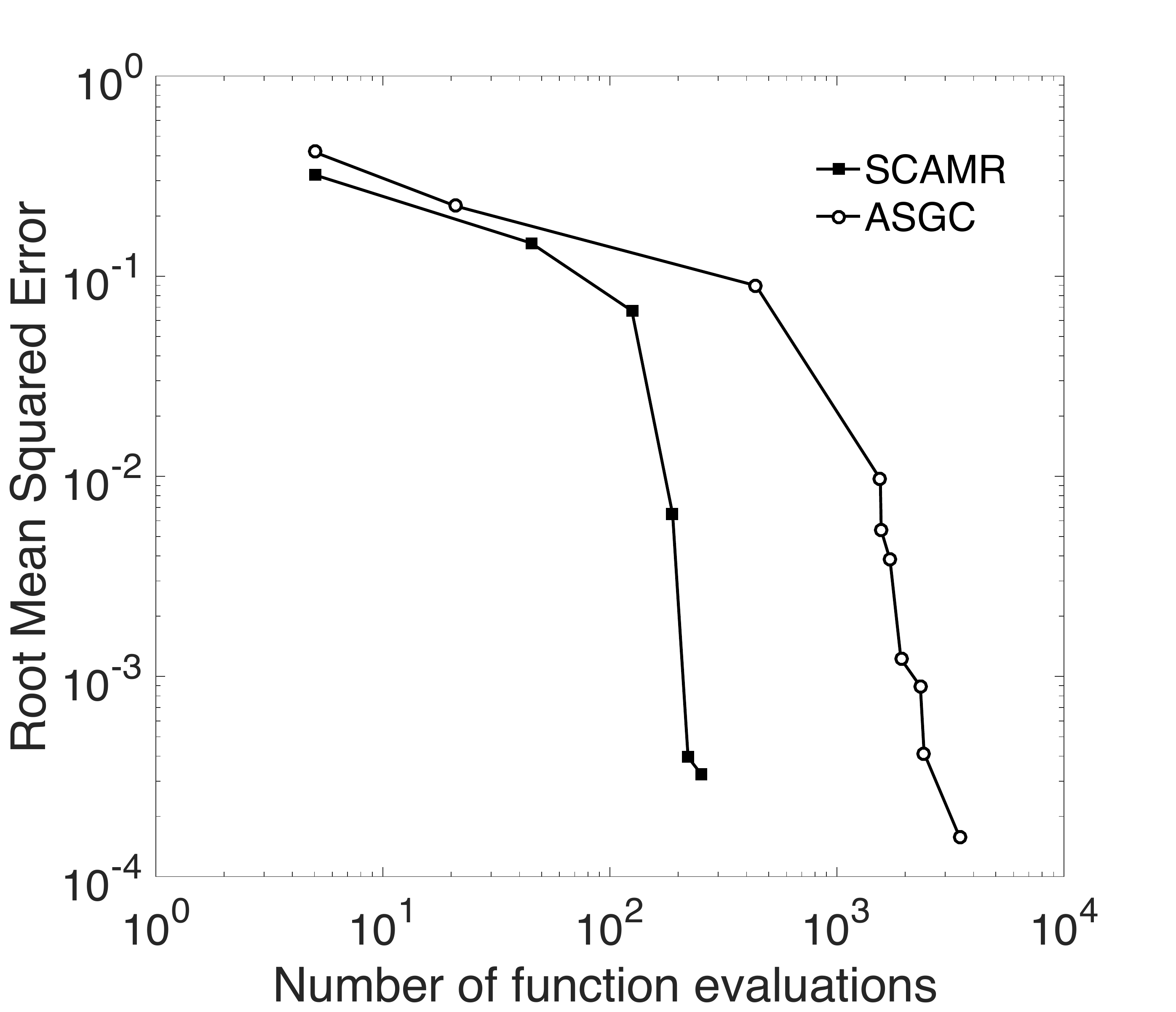} 
\caption{} 
\label{fig:subim131}
\end{subfigure}
\hfill
\begin{subfigure}[b]{0.4\textwidth}
\centering
\includegraphics[width=1.2\linewidth, height=4cm]{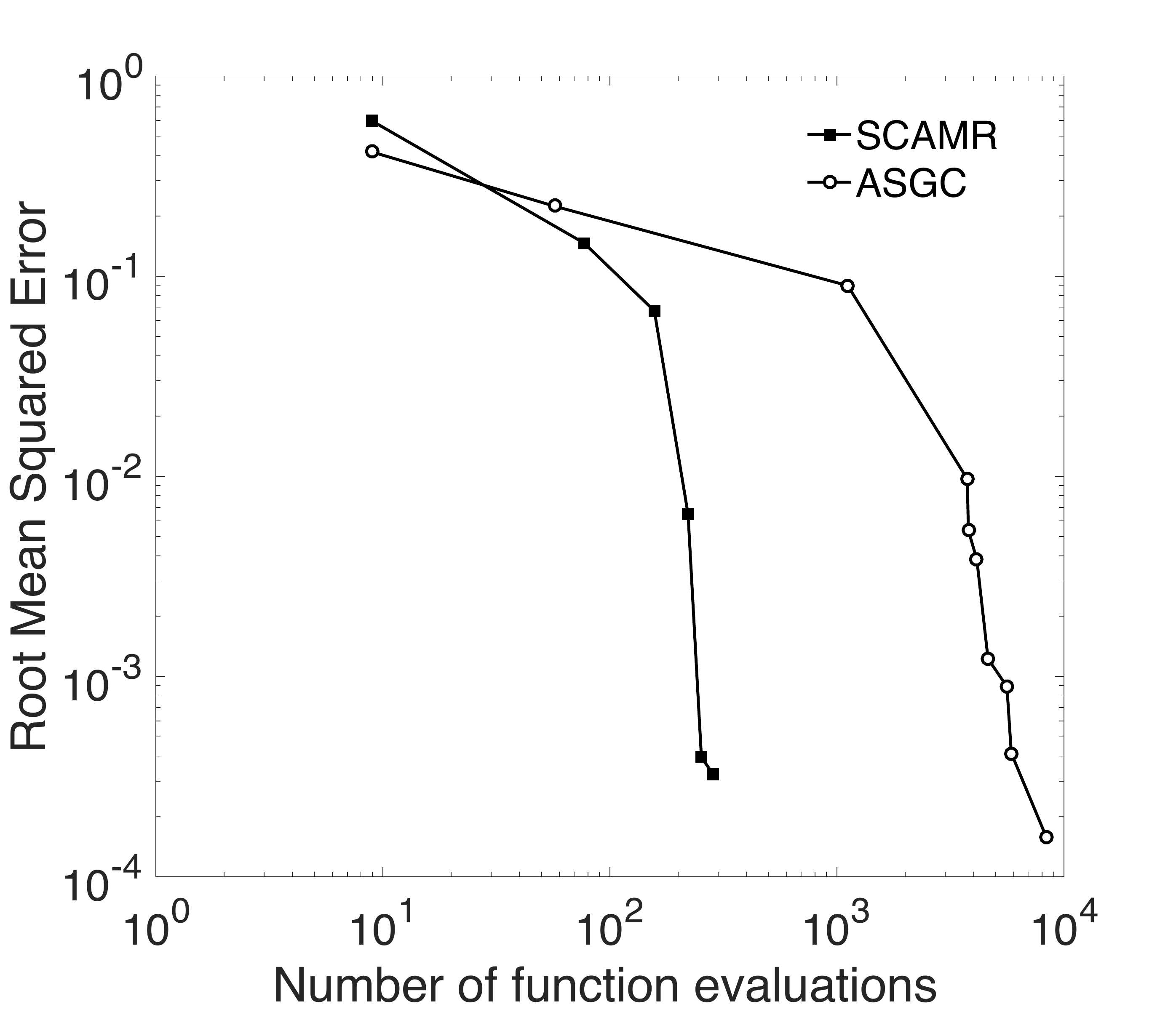} 
\caption{} 
\label{fig:subim131}
\end{subfigure}
\hfill
\begin{subfigure}[b]{0.4\textwidth}
\includegraphics[width=1.2\linewidth, height=4cm]{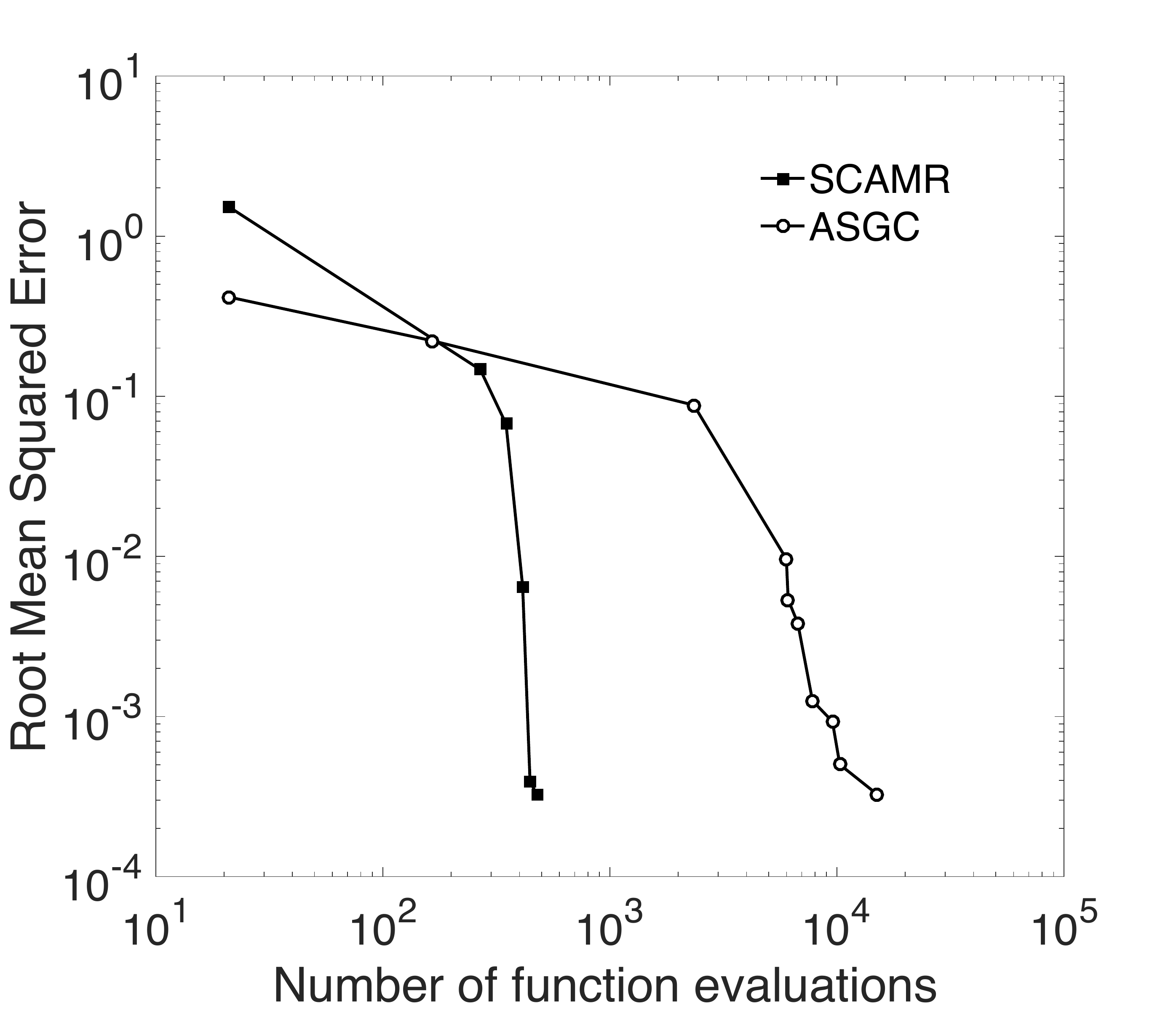} 
\caption{} 
\label{fig:subim131}
\end{subfigure}
\caption{Input domain and error analysis for functions with discontinuty: (a) input domain for function $f_{10}$, (b) numerical error for 2D $f_{10}$, (c) numerical error for 4D $f_{11}$, and (d) numerical error for 10D $f_{12}$. }
\label{fig:imagedisc}
\end{figure}

The proposed SCAMR approach is implemented for these 2-D, 4-D and 10-D functions. The function evaluation locations for 2-D function $f_{10}$ are plotted in Fig. \ref{fig:imagedisc}a, and the error analysis of the numerical approximation from SCAMR for $f_{10}$, $f_{11}$ and $f_{12}$ are provided in Fig. \ref{fig:imagedisc}(b-d).  The numerical approximations are compared to those from ASGC method. From the results, similar conclusions to the previous example can be drawn. 

\subsubsection{SCAMR in a Stochastic Elliptic Problem}
Finally, we apply the SCAMR approach to a stochastic elliptic problem as in \cite{nobile2008sparse,ma2009adaptive}. The model problem is given as    
\begin{align}
-\triangledown (a_n(\omega,x) \triangledown u (\omega,x,y))=f(x,y), \text{in} \ D \times \Gamma \nonumber \\
u(\omega,x,y)=0, \ \text{on}  \ \partial D \times \Gamma
\end{align}
where spatial variable $(x,y) \in D =[0,1]^2$, random variable $\omega \in \Gamma$,  $f(x,y)=\cos(x)\sin(y)$.\\
The diffusion coefficient $a_n(\omega,x)$ is assumed to be a random field that can be approximated in a finite $n$-dimensional stochastic space as:
\begin{equation}
\log(a_n(\omega,x)-0.5)=1+Y_1(\omega)(\frac{\sqrt{\pi}L}{2})^{1/2}+\sum_{i=2}^{n} \xi_i \phi_i(x) Y_i(\omega),
\label{Eq: diff}
\end{equation}
where $Y_i(\omega)$ $[i=1,2, \hdots, n]$ are independent random variables which are uniformly distributed in $[-\sqrt{3},\sqrt{3}]$, and
\begin{equation}
\xi_i=(\sqrt{\pi}L)^{1/2} \exp(\frac{-( \lfloor {\frac{i}{2}}\rfloor  \pi L)^2}{8}), \text{ if } i>1
\end{equation}
and\\
\[
    \phi_i(x):= 
\begin{cases}
   \sin(\frac{\lfloor {\frac{i}{2}}\rfloor  \pi x}{L_p}),			& \text{if i is even},\\
    \cos(\frac{\lfloor {\frac{i}{2}}\rfloor  \pi x}{L_p}),              & \text{if i is odd}
\end{cases}
\]
where $L_p=\max\{1,2L_c\}$, and $L=\frac{L_c}{L_p}$ where $L_c=0.5$ is the correlation length. 

\begin{figure}
\begin{subfigure}[b]{0.4\textwidth}
\centering
\includegraphics[width=1\linewidth, height=4cm]{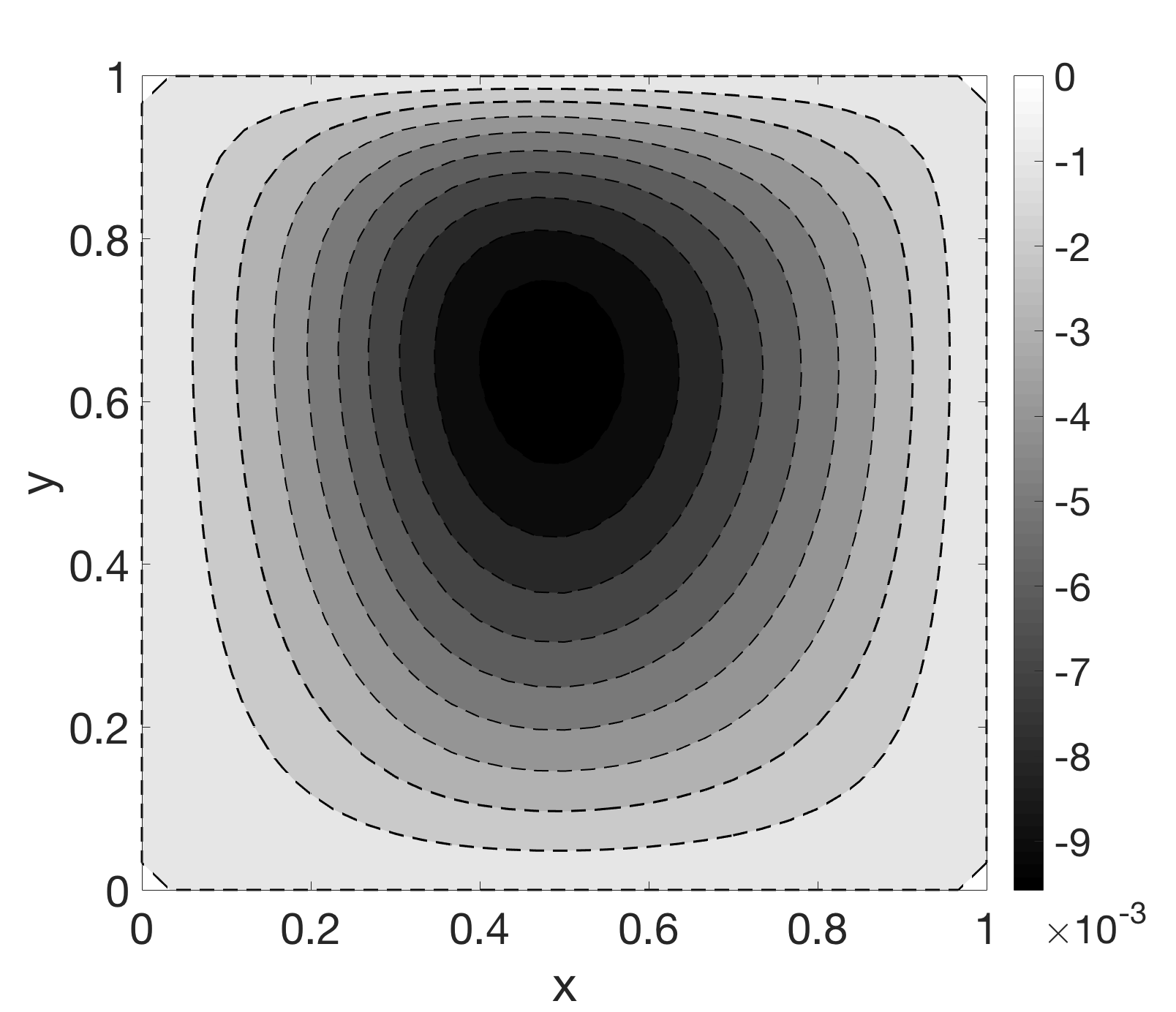} 
\caption*{} 
\label{fig:subim82}
\end{subfigure}
\hfill 
\begin{subfigure}[b]{0.4\textwidth}
\centering
\includegraphics[width=1\linewidth, height=4cm]{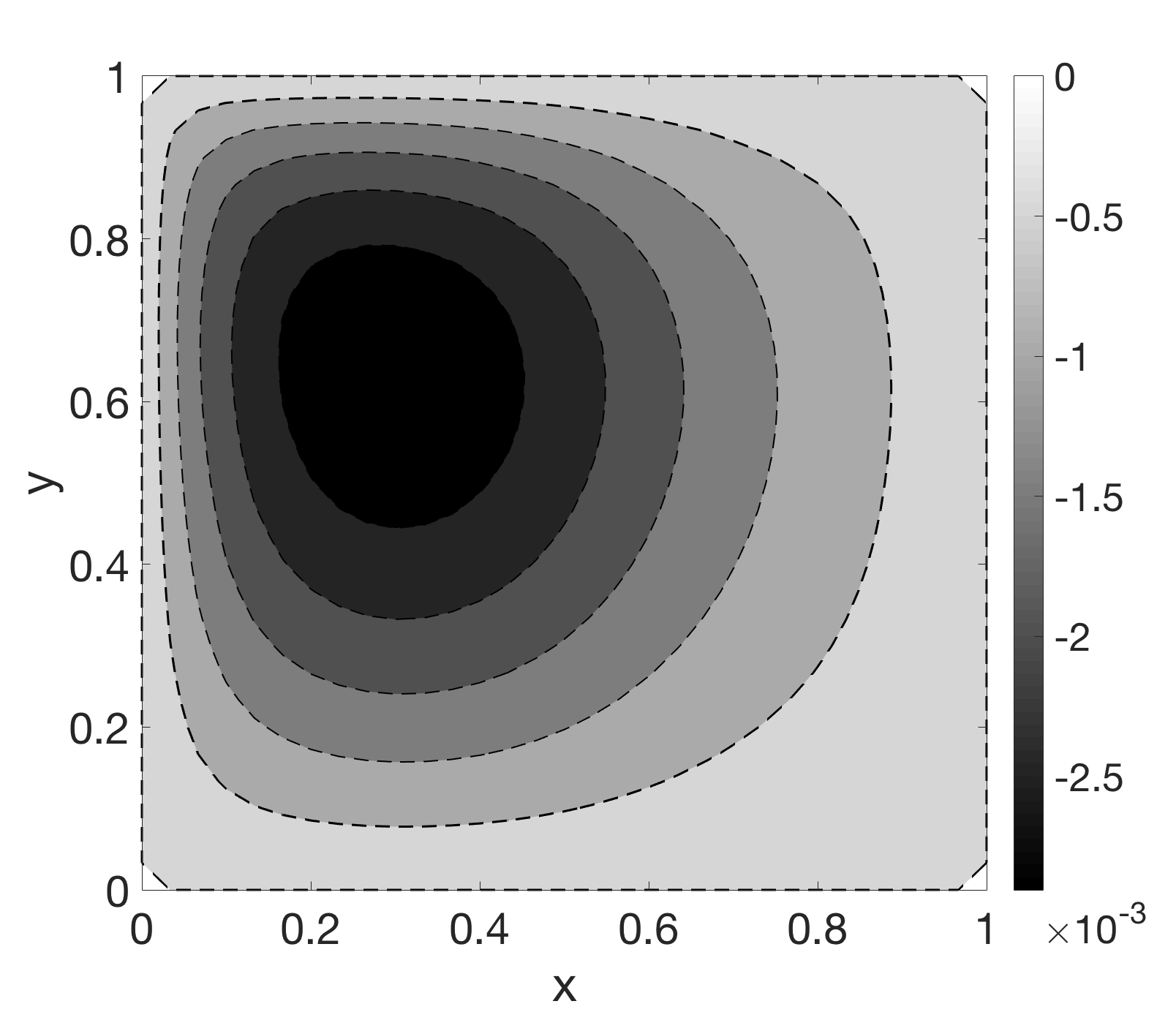} 
\caption*{} 
\label{fig:subim83}
\end{subfigure}
\caption{Two realizations of the output $u$ for  $n=50$ and correlation length $L_c = 0.5$.}
\label{fig:image8}
\end{figure}

\begin{figure}
\centering
\begin{subfigure}[b]{0.4\textwidth}
\centering
\includegraphics[width=1.2\linewidth, height=4cm]{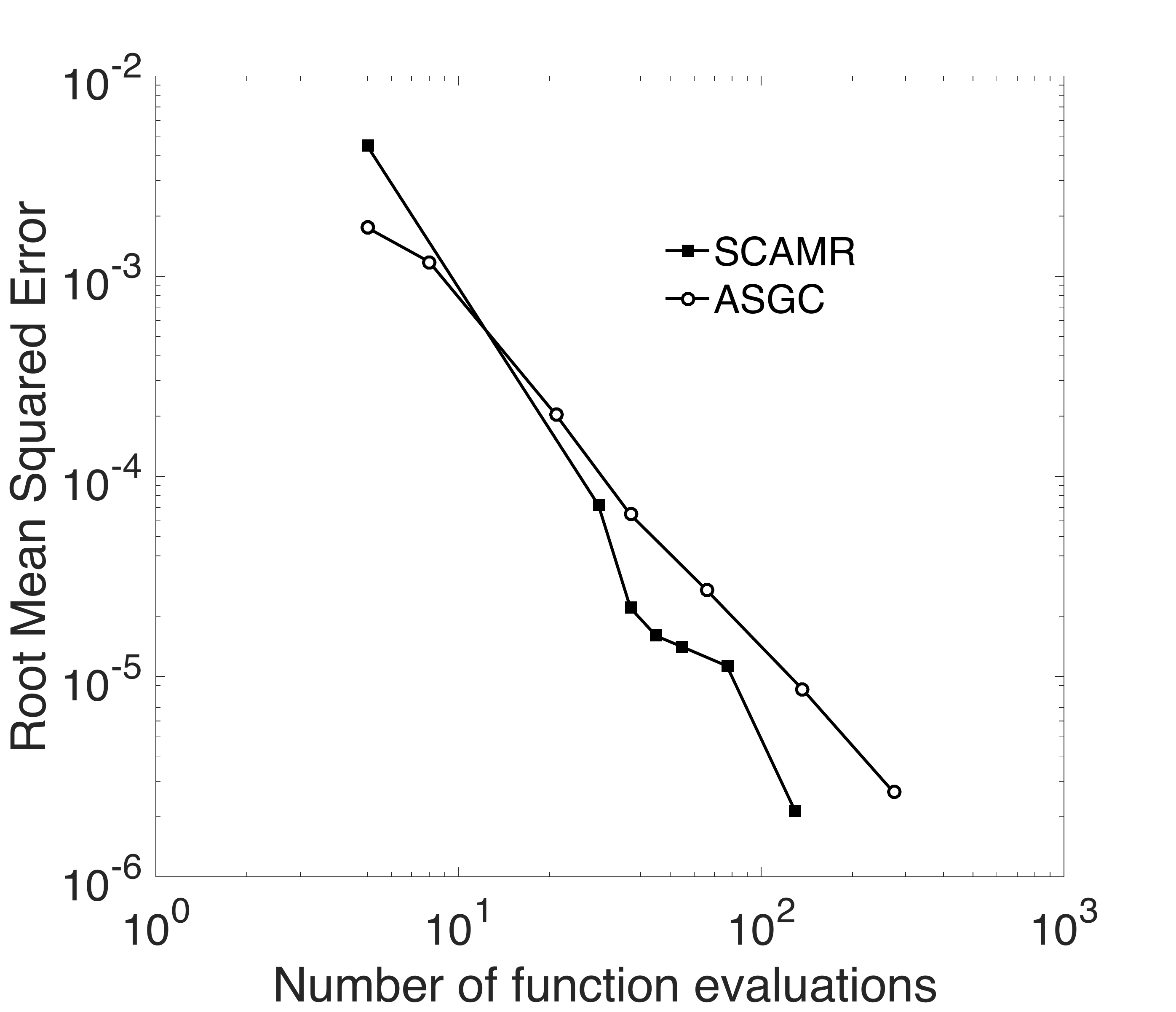} 
\caption{} 
\label{fig:subim131}
\end{subfigure}
\hfill 
\begin{subfigure}[b]{0.4\textwidth}
\centering
\includegraphics[width=1.2\linewidth, height=4cm]{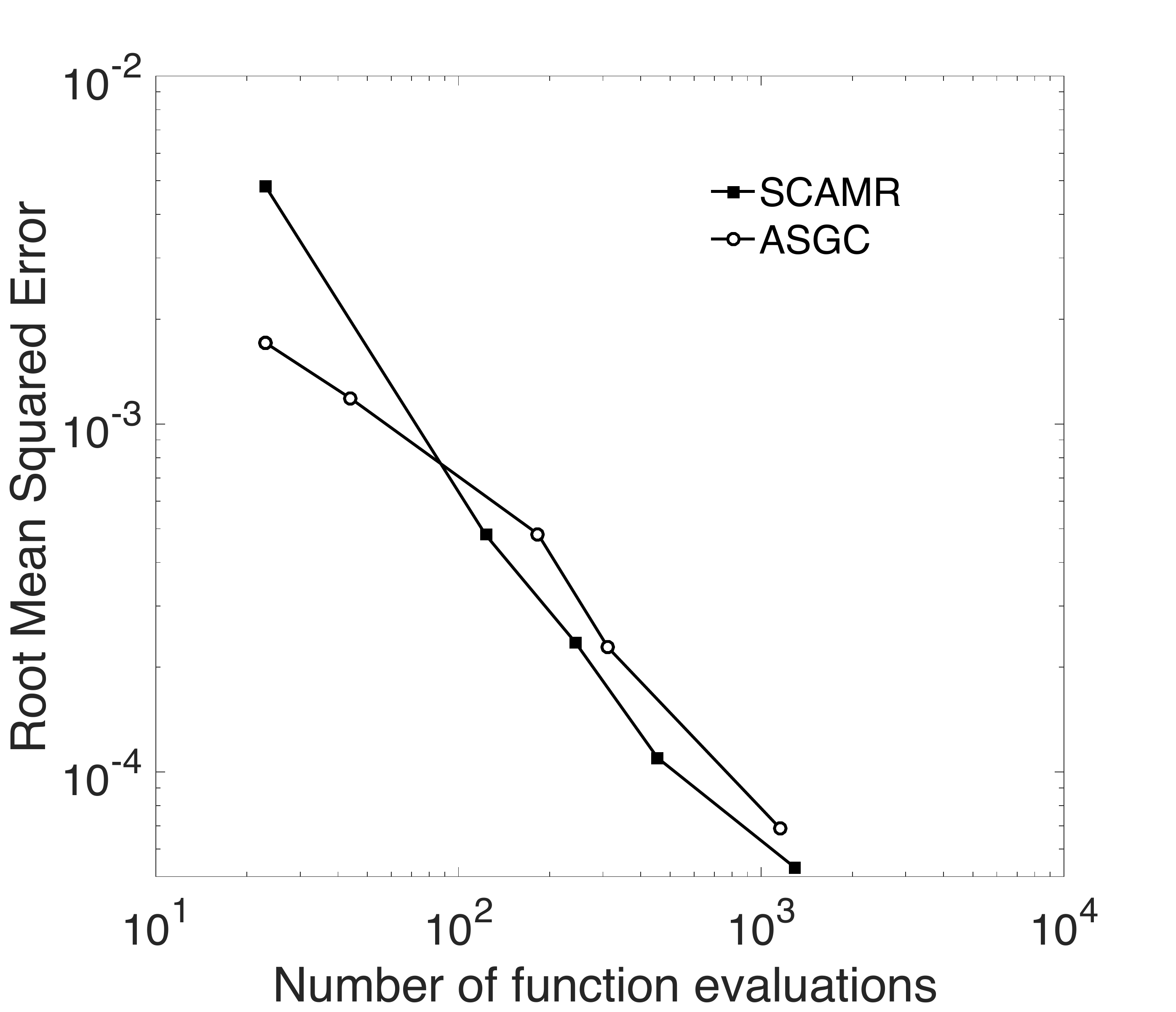} 
\caption{} 
\label{fig:subim131}
\end{subfigure}
\hfill 
\begin{subfigure}[b]{0.4\textwidth}
\centering
\includegraphics[width=1.2\linewidth, height=4cm]{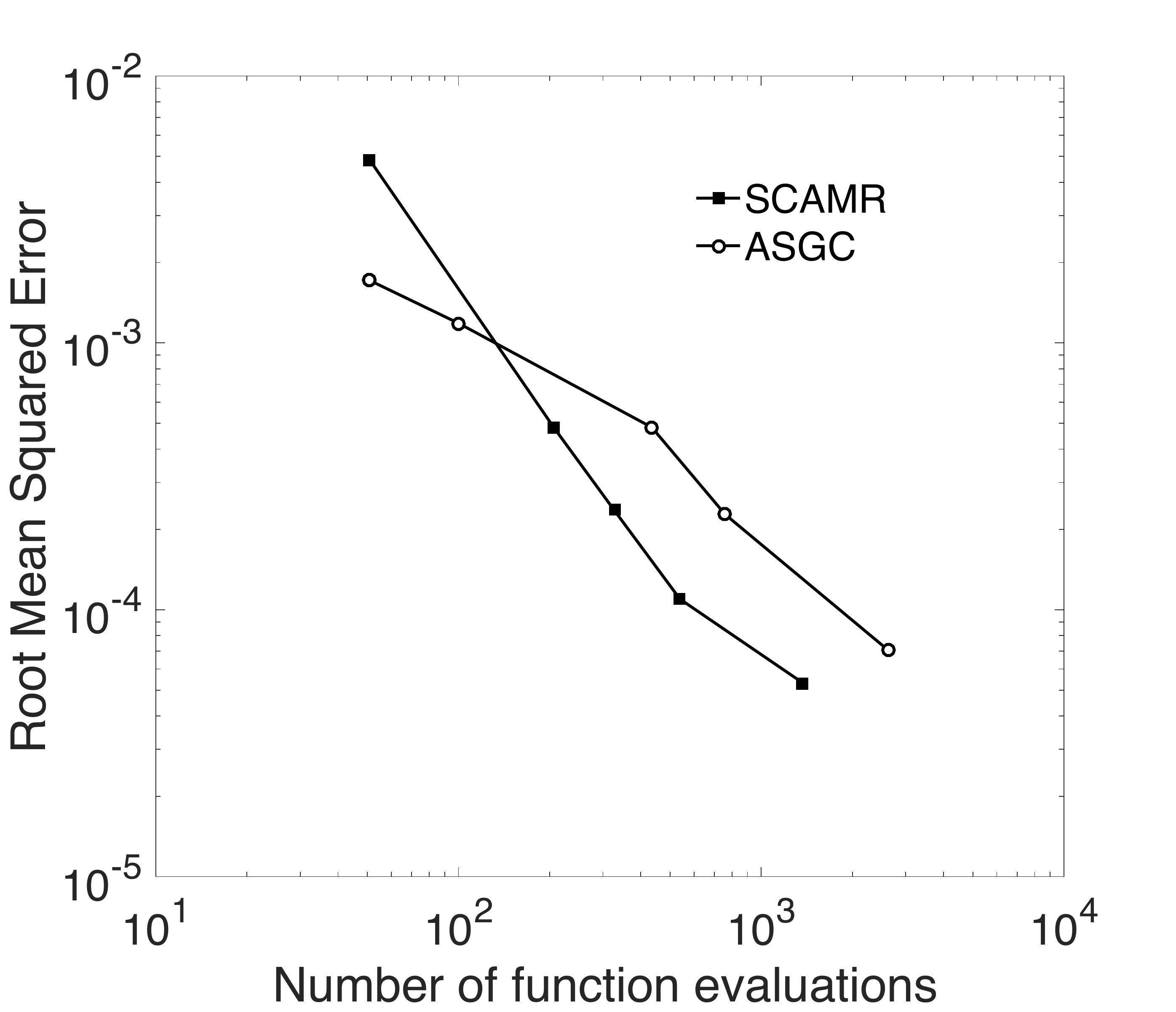} 
\caption{} 
\label{fig:subim131}
\end{subfigure}
\hfill 
\begin{subfigure}[b]{0.4\textwidth}
\centering
\includegraphics[width=1.2\linewidth, height=4cm]{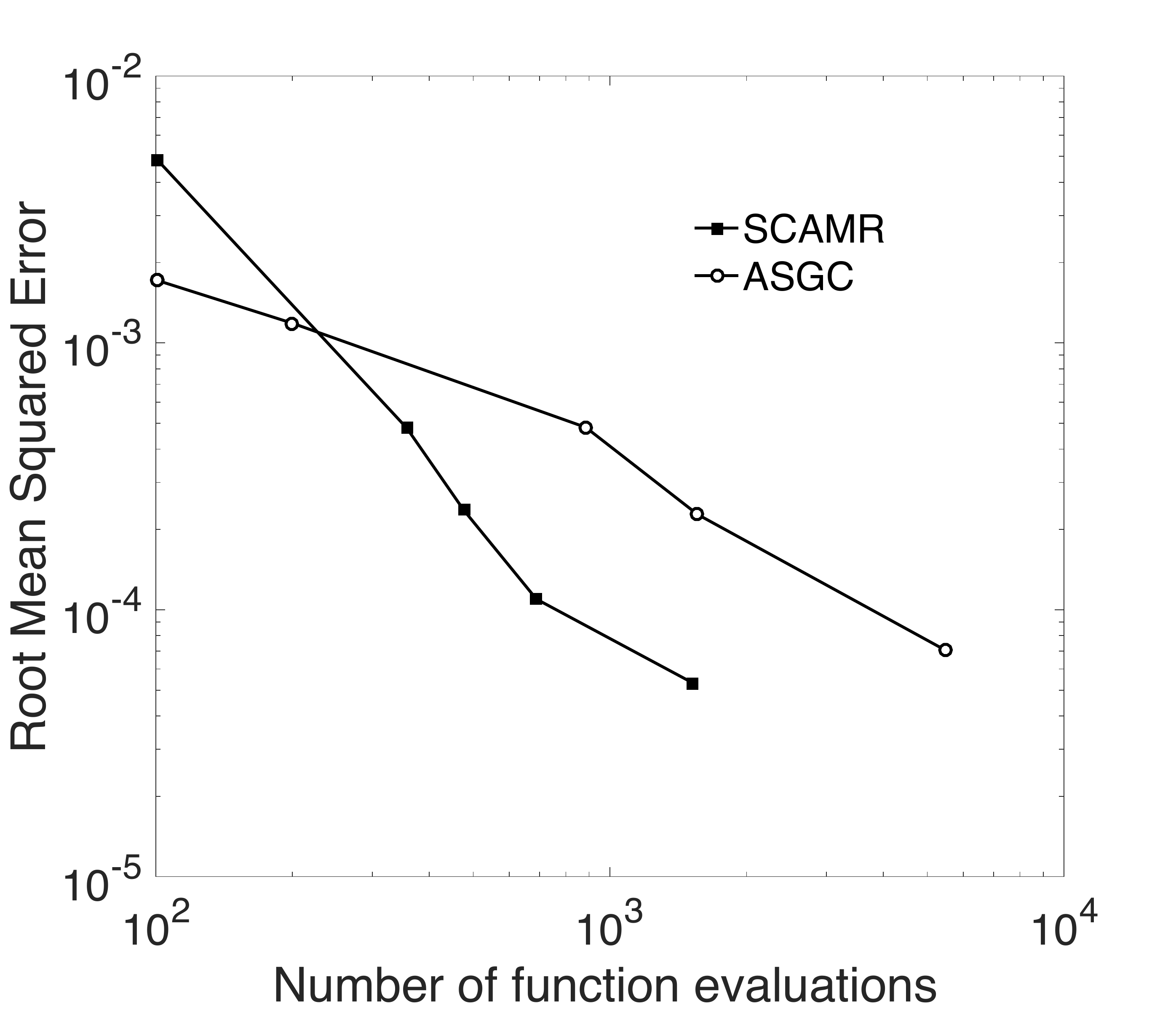} 
\caption{} 
\label{fig:subim131}
\end{subfigure}
\hfill 
\begin{subfigure}[b]{0.4\textwidth}
\centering
\includegraphics[width=1.2\linewidth, height=4cm]{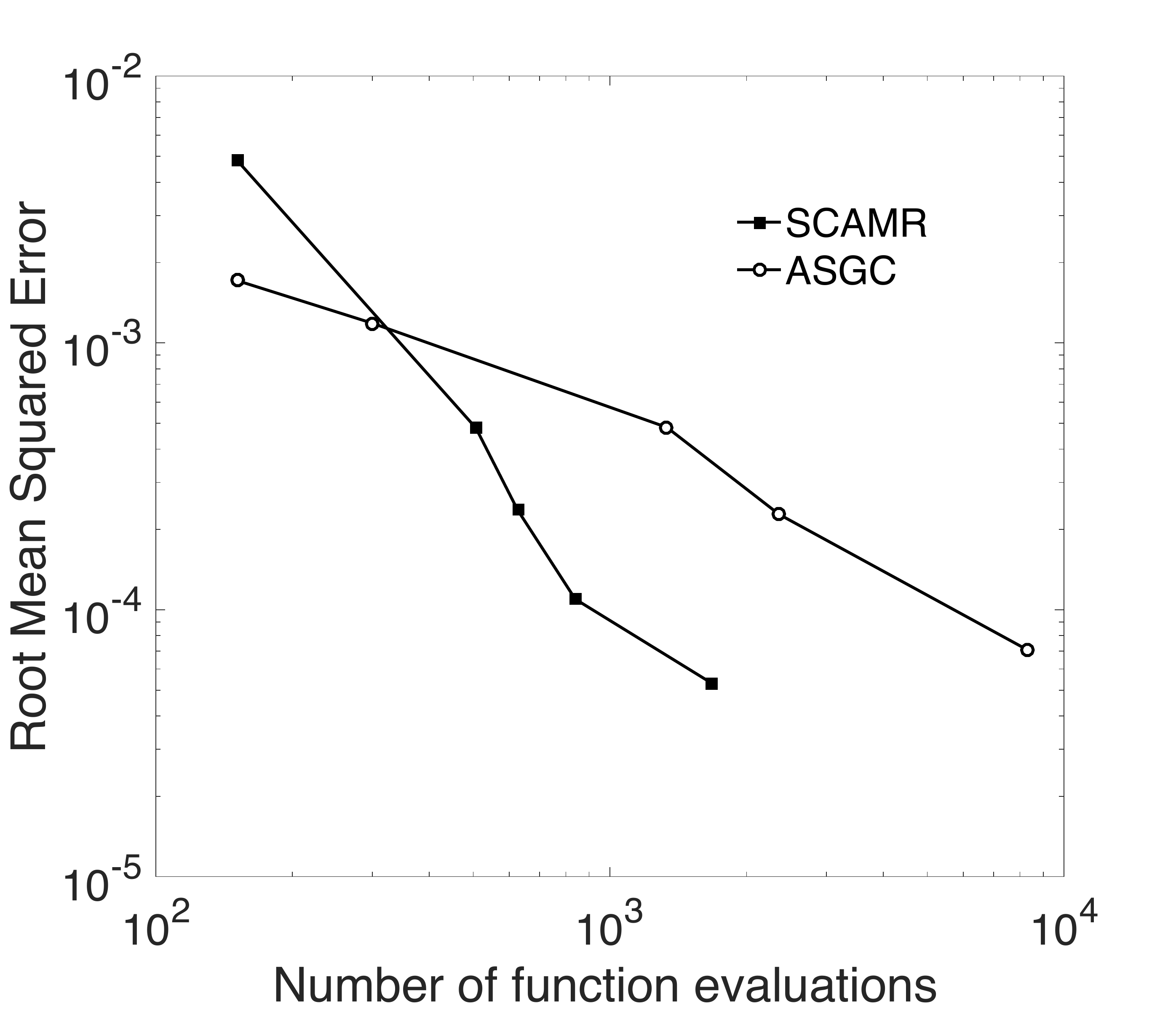} 
\caption{} 
\label{fig:subim131}
\end{subfigure}
\caption{Error analysis of the stochastic elliptic problem with (a) $n = 2$, (b) $n=11$, (c) $n=25$, (d) $n=50$, (e) $n=75$ dimensions for correlation length $L_c = 0.5$.}
\label{fig:image9}
\end{figure}

Without loss of generality, we consider the uncertainty in the output at a fixed point in space $x=y=0.5$, which is the center of the spatial domain. 
Figure \ref{fig:image8} displays two realizations of the output contour in the spatial domain for $n=50$ using the deterministic code of the elliptic problem. The proposed SCAMR approach is implemented for the stochastic elliptic problem with different dimensions $n$ in the random space. The error analysis of the numerical approximations are provided in Fig. \ref{fig:image9}(a-e) for $n=2,11,25, 50, 75$ respectively. The numerical approximations from SCAMR are compared to those from the ASGC method.  From the figure, one can observe that the numerical approximation from SCAMR converges faster for very low dimension such as $n=2$, but it achieves similar convergence rates for large dimensions such as $n=25, 50, 75$. The reason is that the tail terms of Eq. \ref{Eq: diff} for $n>25$ could be negligible due to the fast decay of the eigenvalues $\xi_i$. As with the previous examples, SCAMR converges faster or at a similar rate as ASGC for this problem. 
\subsection{Comparison to HDMR Guided Algorithms for High Dimensional Problems}
To further illustrate the efficiency of SCAMR regarding the model reduction criterion, we implement our proposed approach for more high-dimensional problems and compare the results to those from HDMR-ASGC and MEPCM-A methods. 

A 10-dimensional function is considered to compare the efficiency of SCAMR and HDMR-ASGC \cite{ma2010adaptive}. The error estimate used here is the normalized $L_2$ interpolation error given by
\begin{equation}
\epsilon=\frac{\sqrt{\sum_{i=1}^N (f(\mathbf{x_i})-\tilde{f}(\mathbf{x_i}))^2}}{{\sqrt{\sum_{i=1}^N f(\mathbf{x_i})^2}}},
 \end{equation}
where $f$ is the exact function, $\tilde{f}$ is the numerical approximation using HDMR-ASGC or SCAMR and $N=10^6$ randomly generated samples.

A high dimensional integration problem is then used as an example to compare the SCAMR and the MEPCM-A methods. The error estimate used here is the mean relative error \cite{foo2010multi} given by
\begin{equation}
\epsilon=\frac{|I_{exact}-I_{approx}|}{I_{exact}}
 \end{equation} 
 where $I_{exact}$ is the true mean of the problem and $I_{approx}$ is the numerical approximation of the mean using either MEPCM-A or SCAMR.
\subsubsection{Comparison to HDMR-ASGC}
We consider a 10-D function \cite{ma2010adaptive} given by
\begin{eqnarray}
f_{13}(\bm{x})&=\frac{1}{1+\sum_{i=1}^{10}\alpha_i x_i}
\end{eqnarray}
where parameters $\alpha_i={0.1}/{2^{i-1}}$,  random input $x_i=\sigma y_i$ and $y_i$ are i.i.d. uniform random variables in $\big[-\sqrt{3},\sqrt{3} \ \big]$, $i \in \{1,2,\dots,10\}$. Parameter $\sigma$ is related to the standard deviation of the input and for this example, $\sigma=2$. The weights drop drastically with increase in dimensions and hence the number of effective dimensions is low compared to 10 nominal dimensions. Table \ref{tab:table1} shows a comparison of the normalized $L_2$ interpolation error and the number of points needed for the HDMR-ASGC and the SCAMR approach. It can be seen from the results that SCAMR proves to be more efficient than HDMR-ASGC in approximating $f_{13}$. The HDMR-ASGC results are read directly from Fig. 8 (right) in \cite{ma2010adaptive}. Identification of the low effective dimensions using the interaction check in the SCAMR approach is achieved at a lower computational cost compared to the corresponding check in HDMR-ASGC \cite{ma2010adaptive}. The subsequent surrogate construction of the sub-dimensional problems also requires lesser number of samples when using the second order gPC approximation in SCAMR compared to the linear basis interpolation in the HDMR-ASGC approach. For example, the number of points needed for an $L_2$ error of approximately $6\times10^{-5}$ is around $1575$ points in the case of HDMR-ASGC while the number of points needed for an $L_2$ error of $2.2921 \times 10^{-5}$ using SCAMR is $407$.

\begin{table}[h!]
  \centering
  \caption{HDMR-ASGC and SCAMR error and cost comparison for function $f_{13}$}
  \label{tab:table1}
  \begin{tabular}{cccc}
    HDMR-ASGC &  & SCAMR & \\
    \toprule
    $L_2$ error & Number of points & $L_2$ Error & Number of points\\
    \midrule
      $ \approx 9\times10^{-3}$ & $\approx 200$ &  & \\
      $\approx 1\times10^{-3}$  & $\approx 700$ & $3.9163\times10^{-4}$ & $101$\\
      $\approx 1\times10^{-4}$  & $\approx 1144$ & $8.3553\times10^{-5}$ & $165$\\
     $\approx 6\times10^{-5}$  & $\approx 1575$ & $2.2921\times10^{-5}$ & $407$\\
    \bottomrule
  \end{tabular}
\end{table}
\subsubsection{Comparison to MEPCM-A}
We consider a discontinuous Genz function given by:
\[
    f_{14}(\textbf{x})= 
\begin{cases}
    0,& \text{if } x_1 \geq 0.5 \text{ or } x_2 \ge 0.5,\\
    \exp(\sum_{i=1}^n c_ix_i),              & \text{otherwise}
\end{cases}
\]
where $c_i=e^{-35i/(n-1)}$. Using SCAMR, we evaluate the high dimensional integration $I_{approx}=\int \tilde{f}_{14} (\textbf{x}) d \textbf{x}$ where $\tilde{f}_{14}$ is the numerical approximation to $f_{14}$. The relative mean error is then calculated and compared with MEPCM-A results in Table 6 given in \cite{foo2010multi} with different dimensions $n=100$, $200$, and $300$.
\begin{table}[h!]
  \centering
  \caption{MEPCM-A and SCAMR error and cost comparison for function $f_{14}$}
  \label{tab:table2}
  \begin{tabular}{ccccc}
    & MEPCM-A &  & SCAMR & \\
    \toprule
    $n$ & Relative error & Number of points & Relative Error & Number of points\\
    \midrule
    $100$ & $O(1)$ & $103$ & $2.1308$ & $201$\\
     & $0.0197$  & $20,801$ & $0.0026$ & $2777$\\
     & $0.0098$  & $4,677,148$ & $0.0017$ & $5909$\\
    $200$ & $O(1)$  & $203$ & $2.3705$ & $401$\\
     & $0.067$  & $81,601$ & $0.0105$ & $5121$\\
     & $0.047$  & $36,714,298$ & $0.0021$ & $8397$\\
     $300$ & $O(1)$  & $303$ & $ 2.5435 $ & $601$\\
     & $0.12$  & $182,401$ & $0.7468$ & $2507$\\
     & $0.09$  & $123,111,448$ & $0.01985$ & $21904$\\
    \bottomrule
  \end{tabular}
\end{table}
It can be seen from the form of function $f_{14}$ that the importance of the dimensions decrease exponentially with increase in dimensions. Thus this is an example where the function has a high nominal dimensionality but low effective dimensionality depending on the error tolerance. Table \ref{tab:table2} shows a comparison of the mean relative error and the number of points needed for the MEPCM-A approach and the SCAMR approach. For the SCAMR approach, mean value extraction is performed by generating weighted Clenshaw-Curtis sparse grid points in each of the elements in each subproblem. Then local means are calculated for each subproblem by assigning weights to each element according to their hypervolume. Local means are finally combined together to get the global mean. It can be seen from the results that SCAMR proves to be very efficient in identifying the low effective dimensions. In MEPCM-A, the effective dimensions depend on the parameter $\nu$. Even though $\nu$ is chosen to be small ($\nu=2 \ or \ 3$), the number of terms in HDMR becomes very large for high nominal dimensions. SCAMR thus achieves much better precision with less number of points compared to the MEPCM-A approach. For example, for the $300$-dimensional case, the number of points needed for a relative error of $0.09$ is around $123$ million points in the case of MEPCM-A while the number of points needed for a relative mean error of around $0.02$ is around $22,000$.

\section{Conclusion}

In this paper, an efficient stochastic collocation method with adaptive mesh refinement has been proposed. Specifically, this approach utilizes generalized polynomial chaos as the basis and solves the gPC coefficient using the least squares method, which provides more flexibility on the number and locations of function evaluations. It also implements adaptive mesh refinement to track the discontinuities, and the adaptive criteria of the mesh refinement to check for abrupt variations in the output based on error measured from a second order gPC. In addition, this approach uses a criterion to check possible dimensionality reduction and decomposes the full-dimensional problem to a number of lower-dimensional subproblems, based on the HDMR method. Therefore, for a specific problem, the highest dimensionality of subproblems which involve interacting dimensions, are automatically provided. The effectiveness of this method has been shown using different low and high dimensional, smooth and non-smooth examples. It is noticeable that this approach is particularly efficient for high nominal dimensional problems, like the stochastic elliptic problem with a large number of terms for the diffusivity coefficient, where a significant number of dimensions can be less important (low effective dimensions) and thus non-interacting with other more important dimensions. However, if the dimensions are all coupled in their contribution towards the output of interest, then the efficiency of this method decreases with the increase in the dimensionality of the problem, especially when the response surface is highly non-linear. This is because of the generation of a large number of high dimensional subdomains, where new input points are to be generated according to the sparse grid quadrature. When there is significant non-linearity, the subdomains generally do not converge with the low-order gPC approximation and hence split up into further smaller domains.


\section*{References}
\bibliographystyle{ieeetr}
\bibliography{reference}

\end{document}